\documentclass[a4paper,leqno]{article}
\usepackage[frenchb,english]{babel}
\usepackage{graphicx}
\usepackage{version}

\newcommand{\noi}{\noindent }
\newcommand{\noib}{\noindent $\bullet $~}
\newcommand{\pointir}{\discretionary{.}{}{.\kern.35em---\kern.7em}}
\newcommand{\ie}{\emph{i.e. }}

\newcommand{\resp}{\emph{resp. }}

\newcommand{\como}[1]{\bigskip $\blacktriangleright$ \textbf{ \texttt{[#1]}}
\bigskip}

\newcommand{\comf}[1]{\bigskip \textbf{ \texttt{[#1]}} $\blacktriangleleft$ \bigskip}

\usepackage{latexsym}
\usepackage{amsmath}
\usepackage{amssymb}

\newcommand{\Mh}{\widehat{M}}

\newcommand{\gh}{\widehat{g}}

\newcommand{\Kh}{\widehat{K}}
\newcommand{\Rh}{\widehat{R}}

\newcommand{\rich}{\widehat{\mathrm{Ric}}}

\newcommand{\ind}{\mathrm{Ind}}
\def\buildo#1\over#2{\mathrel{\mathop{\null#2}\limits^{#1}}}
\def\buildu#1\under#2{\mathrel{\mathop{\null#2}\limits_{#1}}}

\newcommand{\B}{\mathbb{B}}

\newcommand{\R}{\mathbb{R}}

\newcommand{\HH}{\mathbb{H}}
\newcommand{\PP}{\mathbb{P}}
\newcommand{\Ss}{\mathbb{S}}

\newcommand{\V}{\mathbb{V}}

\newcommand{\cC}{{\mathcal C}}
\newcommand{\cD}{{\mathcal D}}

\newcommand{\cK}{{\mathcal K}}

\newcommand{\cM}{{\mathcal M}}

\newcommand{\bra}{\langle}
\newcommand{\ket}{\rangle}

\newcommand{\ens}[2]{\{ #1 ~|~ #2 \}}

\newcommand{\pf}{\par{\noindent\textbf{Proof.~}}}

\newbox\qedbox
\setbox\qedbox=\hbox{$\Box$}
\newcommand{\qed}{\hfill\penalty10000\copy\qedbox\par\medskip}
\newtheorem{thm}{Theorem}[section]
\newtheorem{thm-def}[thm]{Theorem and Definition}

\newtheorem{prop}[thm]{Proposition}
\newtheorem{proper}[thm]{Properties}
\newtheorem{prop-def}[thm]{Proposition and Definition}
\newtheorem{lem}[thm]{Lemma}

\newtheorem{cor}[thm]{Corollary}

\title{Minimal hypersurfaces in $\HH^n \times \R$,\\ total curvature and index}

\author{Pierre B\'{e}rard and Ricardo Sa Earp}
\date{November 2009}

\includeversion{pb1-figs}

\excludeversion{pb2}

\excludeversion{pb3}

\excludeversion{pb4}

\begin{document}
\maketitle

\thispagestyle{empty}

\vspace{2cm}

\begin{abstract}
\noi In this paper, we consider minimal hypersurfaces in the product
space $\mathbb{H}^n \times \mathbb{R}$. We begin by studying
examples of rotation hypersurfaces and hypersurfaces invariant under
hyperbolic translations. We then consider minimal hypersurfaces with
finite total curvature. This assumption implies that the
corresponding curvature goes to zero uniformly at infinity. We show
that surfaces with finite total intrinsic curvature have finite
index. The converse statement is not true as shown by our examples
which also serve as useful barriers.
\end{abstract}\bigskip

\bigskip


\newpage

\section{Introduction}\label{S-intro}
In this paper, we focus on complete oriented minimal hypersurfaces
$M$ immersed in $\HH^n\times \R$ equiped with the product metric.
\bigskip

In Section \ref{S-examp}, we study the family $\{\cC_a, a > 0$\} of
hypersurfaces invariant under rotations about the vertical geodesic
$\{0\} \times \R \subset \HH^n \times \R$ (``catenoids'') and the
family $\{\cM_d, d>0 \}$ of hypersurfaces invariant under hyperbolic
translations. These examples generalize to higher dimensions some of
the minimal surfaces constructed in \cite{ST05, NSST08, Sa08}.
\bigskip

In particular, we prove that the $n$-dimensional catenoids $\cC_a$
have vertical heights bounded from above by $\pi/(n-1)$
(Proposition~\ref{P-cat3-10a}). In Section~\ref{SS-dim3-catsta}, we
describe the maximal stable rotationally invariant domains on
$\cC_a$ and we prove that the catenoids have index $1$
(Theorem~\ref{T-cat3-20}). We also give an interpretation in terms
of the envelope of the family $\cC_a$ (Corollary~\ref{C-cat3-21}).
Finally, we observe that the half-catenoid $\cC_a \cap \big( \HH^n
\times \R_{+}\big)$ is not maximally stable. \bigskip

We describe the minimal hypersurfaces invariant under hyperbolic
translations in Theorem~\ref{T-trans3-1}. In particular, we find a
hypersurface $\cM_1$ which is a complete non-entire vertical graph
over a half-space bounded by some hyperplane $\Pi$ in $\HH^n \times
\{0\}$. It takes infinite value data on $\Pi$ and zero asymptotic
boundary value data. When $d<1$, the hypersurface $\cM_d$ is an
entire vertical graph. When $d>1$, it is a bi-graph over the
exterior of an equidistant hypersurface of $\HH^n \times \{0\}$.
\bigskip

In Section \ref{S-index}, we consider the relationships between
finiteness of the total curvature and finiteness of the index. In
dimension $2$, we consider the curvature integrals $\int_M |A_M|^2$
and $\int_M |K_M|$, where $A_M$ is the second fundamental form of
the immersion and $K_M$ the Gauss curvature. Finiteness of these
integrals implies that the corresponding curvatures tend to zero
uniformly at infinity; finiteness of the latter implies finiteness
of the index of the Jacobi (stability) operator
(Theorem~\ref{T-dim2}). The converse statements do not hold. On the
one hand, the catenoids $\cC_a$ have finite index although they have
infinite total intrinsic curvature. This is in contrast with the
case of minimal surfaces in Euclidean $3$-space (\cite{FC82}) and
with the case of surfaces with constant mean curvature $1$ in
hyperbolic $3$-space (\cite{CS90, LR98}). Note that catenoids have
finite total extrinsic curvature. On the other hand, the surfaces
invariant under hyperbolic translations are stable graphs, their
curvature goes to zero at infinity although they have infinite total
curvature. The proof we give of Theorem~\ref{T-dim2} relies mainly
on Simons' equation and the de Giorgi-Moser-Nash method which shows
that finite total curvature implies that the curvature tends to zero
uniformly at infinity. We point out that the finiteness of the
intrinsic total curvature has deep consequences. Under this
assumption on $M$, L. Hauswirth and H. Rosenberg (\cite{HR08},
Theorem 3.1) have indeed shown that the total intrinsic curvature is
quantified, that the ends of $M$ are asymptotic to Scherk type
surfaces and obtained a $C^2$-control on the curvature at infinity.
In dimension $n\ge 3$, we give an upper bound of the index in terms
of the total extrinsic curvature (Theorem~\ref{T-dim3}).
\bigskip

In Section \ref{S-appli}, using the catenoids $\cC_a$ as barriers,
we prove some symmetry and characterization results for minimal
hypersurfaces in $\HH^n \times \R$ whose boundary consists of two
congruent convex hypersurfaces in parallel slices
(Theorem~\ref{T-appli-1min}). We point out that the hypersurfaces
$\cM_d$ ($d<1$ and $d=1$) have been used in \cite{ST08, ST09} as
barriers for the Dirichlet problem and that they play a crucial role
for some existence theorem for the vertical minimal surface
equation. \bigskip

Finally, we point out that most of our results may be established if
the ambient space is one of the product spaces $\HH^n \times \R^k$
or $\HH^n \times \HH^k$. \bigskip

The authors would like to thank the Mathematics Department of
PUC-Rio (PB) and the Institut Fourier -- Universit\'{e} Joseph Fourier
(RSA) for their hospitality. They gratefully acknowledge the
financial support of CNPq, FAPERJ, Universit\'{e} Joseph Fourier and
R\'{e}gion Rh\^{o}ne-Alpes.

\section{General framework}\label{S-frame}

\subsection{Notations}\label{SS-not}

We consider hypersurfaces $M$ immersed in the space $\Mh := \HH^n
\times \R$ equiped with the product metric $\gh = g_{\B} + dt^2$,
where $g_{\B}$ is the hyperbolic metric,
\begin{equation}\label{E-rot3-1}
g_{\B} := \big(\frac{2}{1-|x|^2}\big)^2 \big( dx_1^2 + \cdots +
dx_n^2\big).
\end{equation}

We have chosen the ball model $\B$ for the $n$-dimensional
hyperbolic space $\HH^n$.

\subsection{Jacobi operator, Index, Jacobi fields} \label{SS-jac}

Let $M^n \looparrowright \Mh^{n+1}$ be an orientable minimal
hypersurface in an oriented Riemannian manifold $\Mh$ with metric
$\gh$. Let $N_M$ be a unit normal field along $M$ and let $A_M$ be
the second fundamental form of the immersion with respect to $N_M$.
Let $\rich$ be the normalized Ricci curvature of $\Mh$. The second
variation of the volume functional gives rise to the \emph{Jacobi
operator} (or stability operator) $J_M$ of $M$ (see \cite{Si68,
Law80, CM02}),
\begin{equation}\label{E-dim2-jac3}
J_M := - \Delta_M - \big( |A_M|^2 + \rich(N_M) \big),
\end{equation}
where $\Delta_M$ is the (non-positive) Laplacian on $M$ (for the
induced metric). \bigskip

Given a relatively compact regular domain $\Omega$ on the
hypersurface $M$, we let $\ind(\Omega)$ denote the number of
negative eigenvalues of $J_M$ for the Dirichlet problem on $\Omega$
(this is well defined because $\Omega$ is compact). The \emph{index}
of $M$ is defined to be the supremum ($\le + \infty$)
\begin{equation}\label{E-dim2-ind}
\ind(M) := \sup \ens{\ind(\Omega)}{\Omega \Subset M},
\end{equation}
taken over all relatively compact regular domains. \bigskip

Let $\lambda_1(\Omega)$ be the least eigenvalue of the operator
$J_M$ with Dirichlet boundary conditions in $\Omega$. Recall that a
relatively compact regular domain $\Omega$ is said to be
\emph{stable}, if $\lambda_1(\Omega) > 0$; \emph{unstable}, if
$\lambda_1(\Omega) < 0$; \emph{stable-unstable}, if
$\lambda_1(\Omega) = 0$. More generally, we say that a domain
$\Omega$ is stable if any relatively compact subdomain is stable.

\begin{proper}\label{P-frame-1}
Recall the following properties.
\begin{enumerate}
    \item Let $\Omega$ be a stable-unstable relatively
    compact domain. Then, any smaller domain is stable while any larger
    domain is unstable (monotonicity of Dirichlet eigenvalues).
    \item Of particular interest are the solutions of the equation
    $J_M(u)=0$. We call such functions \emph{Jacobi fields} on $M$.
    Let $X_a : M^n \looparrowright (\Mh^{n+1},\gh)$ be a one-parameter
    family of oriented minimal immersions, with variation field
    $V_a = \frac{\partial X_a}{\partial a}$ and unit normal $N_a$.
    Then, the function $\gh(V_a,N_a)$ is a Jacobi field on $M$
    (\cite{BGS87}, Theorem 2.7).
    \item Let $\Omega$ be a relatively compact domain on a minimal
    manifold $M$. If there exists a positive function $u$ on $\Omega$
    such that $J_M(u) \ge 0$, then $\Omega$ is stable (\cite{FCS80},
    Theorem 1).
\end{enumerate}
\end{proper}

\section{Examples of minimal hypersurfaces in $\HH^n \times \R$}
\label{S-examp}

In this section we give examples of minimal hypersurfaces in $\HH^n
\times \R$. We use these examples as guidelines and counter-examples
to study the relationships between index properties of the Jacobi
operator and the finiteness of some total curvature of $M$, see
Theorems~\ref{T-dim2} and \ref{T-dim3}. We also use them as barriers
for a symmetry and characterization result in Section~\ref{S-appli}.

\subsection{Rotation hypersurfaces in $\HH^n \times \R$}\label{SS-dim3-rot}

We first consider rotation hypersurfaces about a vertical geodesic
axis in $\HH^n \times \R$. Up to isometry, we can assume the
rotation axis to be $\{0\}\times \R$. Recall that we take the ball
model for $\HH^n$. \bigskip

Take the vertical plane $\V := \ens{(x_1, \ldots , x_n,t) \in
\Mh}{x_1 = \cdots = x_{n-1} = 0}$ and consider a generating curve
$\big( \tanh (f(t)/2), t\big)$ for some positive function $f$ which
represents the hyperbolic distance to the axis $\R$, at height $t$.
\bigskip

We define a \emph{rotation hypersurface} $M \looparrowright \Mh$ by
the ``parametrization''
\begin{equation}\label{E-rot3-3a}
X : \left\{%
\begin{array}{l}
\R_{+} \times S^{n-1} \rightarrow \Mh, \\
(t , \xi ) \mapsto \big( \tanh (f(t)/2)\xi , t\big), \\
\end{array}%
\right.
\end{equation}
where $\xi = (\xi_1, \ldots , \xi_n)$ is a point in the unit sphere
$S^{n-1}$ and $\tanh (\rho/2)\xi$ stands for the point $(\tanh
(\rho/2)\xi_1, \ldots , \tanh (\rho/2)\xi_n)$ in the ball $\B$.
\bigskip

The basic tangent vectors to the immersion $X$ are
$$T(t,\xi ) := T_{t,\xi}X(\partial_{t}) = \big( \frac{f_t(t)}{2 \cosh^2(f(t)/2)}\xi,
1\big),$$ where $f_t$ is the derivative of $f$ with respect to $t$,
and
$$U(t,\xi,u) := T_{t,\xi}X(u) = \big( \tanh (f(t)/2) u,0\big),$$
where $u \in T_{\xi}S^{n-1}$ is a unit vector. \bigskip

\begin{pb3}
\como{PB3}

\input{pb3-04}

\comf{PB3}\bigskip
\end{pb3}

We collect basic formulas in the next proposition whose proof is
straightforward.

\begin{prop}\label{P-rot3-5}
Let $(M, g_M) \looparrowright (\Mh, \gh)$ be an isometric
im\-mersion. We have the following formulas in the parametrization
$X$ on $\R \times S^{n-1}$.
\begin{enumerate}
    \item The induced metric $g_{M}$ is given by
    \begin{equation}\label{E-rot3-5}
    g_{M} = \big( 1 + f_t^2(t) \big) dt^2 +
    \sinh^2(f(t)) g_S \, ,
    \end{equation}
    where $g_S$ is the canonical metric on $S^{n-1}$.
    \item The Riemannian measure $d\mu_M$ for the metric $g_M$ is given by
    \begin{equation}\label{E-rot3-6}
    d\mu_{M} = \big( 1 + f_t^2(t) \big)^{1/2}
    \sinh^{n-1}(f(t)) \, dt \, d\mu_S,
    \end{equation}
    where $d\mu_S$ is the canonical measure on the sphere.
    \item The unit normal field to the immersion can be chosen to be
    \begin{equation}\label{E-rot3-7}
    N_{M}(t,\xi) = (1 + f_t^2(t))^{-1/2} \big(\frac{-1}{2 \cosh^2(f(t)/2)} \xi ,
    f_t(t)\big)\,.
    \end{equation}
    In particular, the vertical component of the unit normal field
    is given by
    \begin{equation}\label{E-rot3-7a}
    v_M(t) := f_t(t)\big(1 + f_t^2(t) \big)^{-1/2}.
    \end{equation}
\end{enumerate}
\end{prop}

At the point $X(t,\xi)$, the principal directions of curvature of
$M$ are
\begin{itemize}
    \item the tangent to the meridian curve in the vertical $2$-plane
    $$\V_{\xi} = \ens{\big( \tanh(\rho/2) \xi , t \big)}{(\rho ,t)
    \in \R^2},$$
    \item the vectors tangent to the distance sphere $X(t, S^{n-1})$
    at $\xi$ in the hyperbolic slice
    $\HH^n \times \{t\}$, where the restriction of the second fundamental
    form $A_M$ is a scalar multiple of the identity.
\end{itemize}
The principal curvatures with respect to $N_M$ are
\begin{itemize}
    \item $k_n(t)$, the principal curvature in the direction
    tangent to the meridian curve, given by
    \begin{equation}\label{E-rot3-8a}
    k_n(t) = - f_{tt}(t) (1 + f_t^2(t))^{-3/2},
    \end{equation}
    \item the principal curvatures in the directions tangent to
    $X(t,S^{n-1})$ at $X(t,\xi)$,
    \begin{equation}\label{E-rot3-8b}
    k_1(t) = \cdots = k_{n-1}(t) = \coth (f(t)) (1 + f_t^2(t))^{-1/2}.
\end{equation}
\end{itemize}

We conclude that the mean curvature $H(t)$ of the rotation
hypersurface $M \looparrowright \Mh$ with respect to the unit normal
$N_M$ is given by
\begin{equation}\label{E-rot3-9a}
n H(t) = - f_{tt}(t) (1 + f_t^2(t))^{-3/2} + (n-1) \coth (f(t)) (1 +
f_t^2(t))^{-1/2},
\end{equation}
or
\begin{equation}\label{E-rot3-9b}
n f_t(t) \, \sinh^{n-1}(f(t))\, H(t) =
\partial_{t} \Big( \sinh^{n-1}(f(t)) (1 + f_t^2(t))^{-1/2}\Big).
\end{equation}

\subsection{Catenoids in $\HH^n \times \R$}\label{SS-dim3-cat}

In this Section, we describe the minimal rotation hypersurfaces
about $\{0\} \times \R$, in $\HH^n \times \R$. By analogy with the
Euclidean case, we call them \emph{catenoids}. They are the higher
dimensional counterparts of the catenoids constructed in
\cite{ST05}.

Given some $a > 0$, let $\big(I_a, f(a,\cdot) \big)$ denote the
maximal solution of the Cauchy problem
\begin{equation}\label{E-cat3-10}
\left\{%
\begin{array}{lll}
f_{tt} & = & (n-1) \coth (f) (1 + f_t^2),\\
f(0) &=& a > 0,\\
f_t(0)&=& 0,\\
\end{array}%
\right.
\end{equation}

where $f_{t}$ and $f_{tt}$ are the first and second derivatives of
$f$ with respect to $t$.\bigskip

\begin{prop}\label{P-cat3-10a}
For $a > 0$, the maximal solution $\big(I_a, f(a,\cdot) \big)$ gives
rise to the generating curve $C_a$, $t \mapsto \big(
\tanh(f(a,t)),t\big)$ (\emph{catenary}), of a complete minimal
rotation hypersurface $\cC_a$ (\emph{catenoid}) in $\HH^n \times
\R$, with the following properties.
\begin{enumerate}
    \item The interval $I_a$ is of the form $I_a = ]-T(a), T(a)[$ for
    some finite positive number $T(a)$ and $f(a,\cdot)$ is an even function
    of the second variable.
    \item For all $t \in I_a$, $f(a,t) \ge a$.
    \item The derivative $f_t(a,\cdot)$ is positive on $]0,T(a)[$,
    negative on $]-T(a),0[$.
    \item The function $f(a,\cdot)$ is a bijection from $[0,T(a)[$
    onto $[a,\infty[$, with inverse function $\lambda(a,\cdot)$ given by
    \begin{equation}\label{E-cat3-10b}
    \lambda(a,\rho) = \sinh^{n-1}(a) \int_a^{\rho} \big( \sinh^{2n-2}(u) -
    \sinh^{2n-2}(a)\big)^{-1/2} \, du.
    \end{equation}
    \item The catenoid $\cC_a$ has finite vertical height $h_R(a)$,
    \begin{equation}\label{E-cat3-10c}
    h_R(a) = 2 \sinh^{n-1}(a) \int_a^{\infty} \big( \sinh^{2n-2}(u) -
    \sinh^{2n-2}(a)\big)^{-1/2} \, du.
    \end{equation}
    \item The function $a \mapsto h_R(a)$ increases from $0$ to
    $\frac{\pi}{(n-1)}$ when $a$ increases from $0$ to infinity.
    Furthermore, given $a \not = b$, the generating catenaries $C_a$ and $C_b$
    intersect at exactly two symmetric points.
\end{enumerate}
\end{prop}

\pf \emph{Assertion 1} follows from the Cauchy-Lipschitz theorem for
some positive $T(a)$ which is finite as we will see below.
\medskip

\emph{Assertion 2} follows from the fact that $\sinh^{n-1}(f(a,t))
\big( 1 + f_t^2(a,t) \big)^{-1/2} = \sinh^{n-1}(a)$ for all $t \in
]-T(a),T(a)[$ (see (\ref{E-rot3-9b})). \medskip

\emph{Assertion 3} is clear. \medskip

\emph{Assertion 4}. According to Assertion 3, $t \mapsto f(a,t)$ is
increasing so that it has a limit when $t$ tends to $T(a)$ and this
limit must be infinite because we took a maximal solution. It
follows that the inverse function $\lambda(a,\cdot)$ maps
$[a,\infty[$ onto $[0,T(a)[$ and that
$\lambda_{\rho}(a,f(a,t))f_t(a,t) \equiv 1$. Finally, we find that
$\lambda_{\rho}(a,\rho)= \sinh^{n-1}(a) \big( \sinh^{2n-2}(\rho) -
\sinh^{2n-2}(a) \big)^{-1/2}$ on $]a,\infty[$ and the formula for
$\lambda(a,\rho)$ follows because $f(a,0)=a$. Note that the integral
(\ref{E-cat3-10b}) converges at $u=a$. \medskip

\emph{Assertion 5}. We have that $h_R(a) = 2 T(a)$, where
$$T(a) = \lim_{\rho \to \infty}\lambda(a,\rho) =
\sinh^{n-1}(a) \int_a^{\infty} \big( \sinh^{2n-2}(u) -
\sinh^{2n-2}(a)\big)^{-1/2} \, du,$$ where the integral converges at
both $a$ and $\infty$.\medskip

\emph{Assertion 6}. By  a change of variables, we can write

$$T(a) = \sinh(a) \int_1^{\infty}\big( v^{2n-2} - 1\big)^{-1/2}
\big( \sinh^2(a) v^2 +1 \big)^{-1/2} \, dv$$

and compute the derivative

$$T'(a) = \cosh(a) \int_1^{\infty}\big( v^{2n-2} - 1\big)^{-1/2}
\big( \sinh^2(a) v^2 +1 \big)^{-3/2} \, dv > 0.$$

Note that

$$\sinh(a) \big( v^{2n-2} - 1\big)^{-1/2}
\big( \sinh^2(a) v^2 +1 \big)^{-1/2} \le v^{-1}\big( v^{2n-2} -
1\big)^{-1/2}$$

and that the right-hand side is in $L^1([1,\infty[)$ for $n\ge 2$,
so that we can take the limits under the integral and obtain that
$\lim_{a \to 0}T(a) = 0$ and $\lim_{a \to \infty}T(a) =
\int_1^{\infty}v^{-1}\big( v^{2n-2} - 1\big)^{-1/2} \, dv$. The last
integral can be calculated explicitly because $\big( \arctan
\sqrt{v^N-1} \big)' = \frac{N}{2 \sqrt{v^N-1}}$. The last assertion
follows by considering the function $\lambda(a,\rho) -
\lambda(b,\rho)$ and by using the monotonicity of $T(a)$. \qed
\bigskip

\begin{figure}[h!tb]
\begin{minipage}[l]{.45\linewidth}
\textbf{Remark}. ~The above proposition shows that the catenoids in
$\HH^n \times \R$ have uniformly bounded finite vertical height.
This is in contrast with the Euclidean catenoids ($n\ge 3$) which
have finite, yet unbounded, vertical heights.
\end{minipage}\hfill %
\begin{minipage}[r]{.40\linewidth}
    \includegraphics[scale=.25, angle=-90]{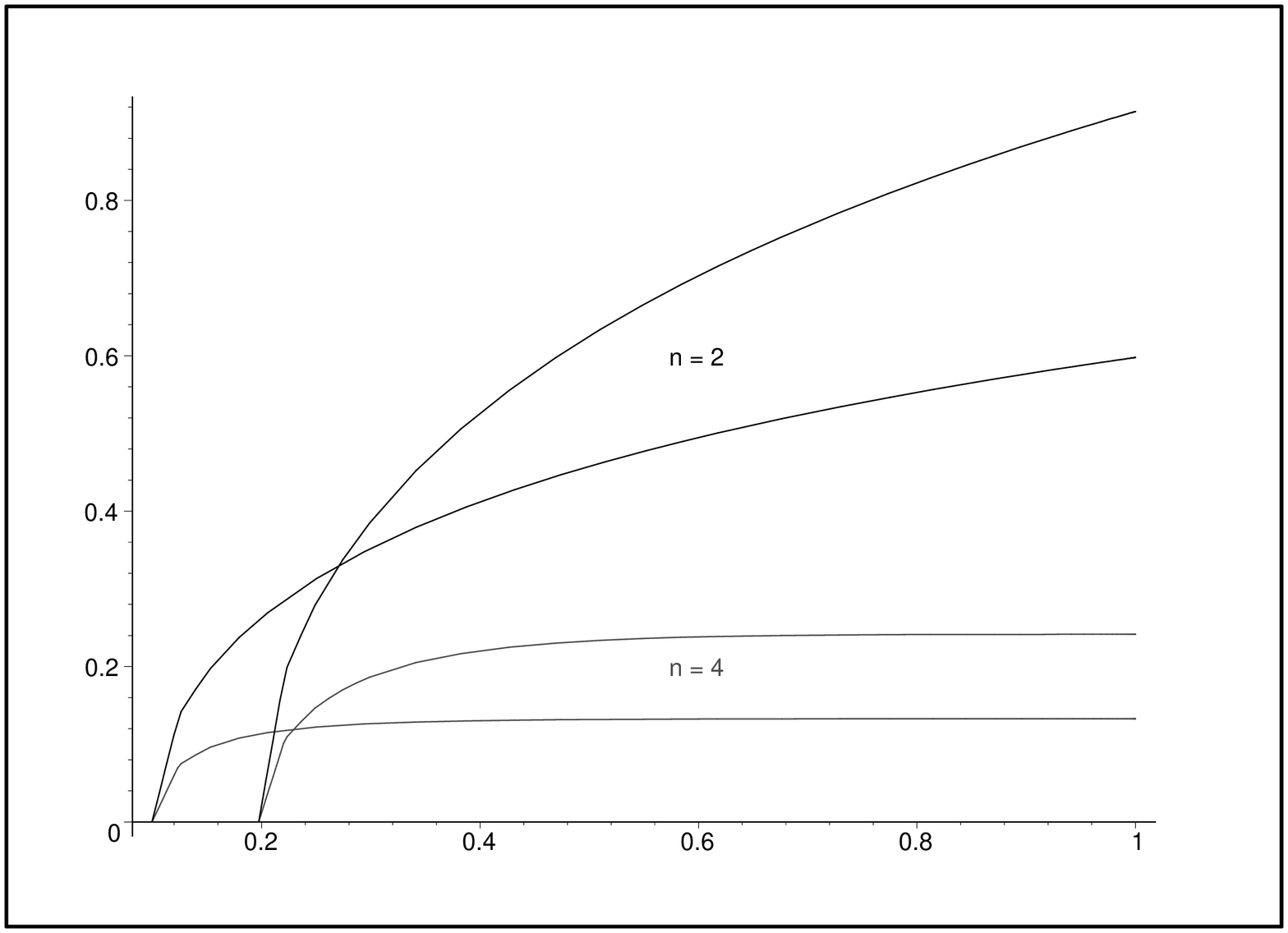}
    \caption[Catenaries $n=2, 4$]{Catenaries $n=2, 4$}
    \label{F-catene4}
\end{minipage}
\end{figure}

\textbf{Notations}. ~Let $N_a$ denote the unit normal to the
catenoid $\cC_a$, let $A_a$ denote its second fundamental form
relative to the normal $N_a$ and let $d\mu_a$ denote its Riemannian
measure. When $n=2$, let $K_a$ denote the Gauss curvature of
$\cC_a$. We state the following proposition for later purposes.

\begin{prop}\label{P-cat3-11a}
For $a>0$, the $n$-dimensional catenoid $\cC_a$ in $\HH^n \times \R$
has infinite volume and finite \emph{total extrinsic curvature}
$\int_{\cC_a} |A_a|^n \, d\mu_a$. When $n=2$, the catenoid $\cC_a$
has infinite \emph{total intrinsic curvature} $\int_{\cC_a} |K_a| \,
d\mu_a$.
\end{prop}

\pf We can restrict to the upper half-catenoid, $\cC_{a,+} = \cC_a
\cap (\HH^n \times \R_{+})$, which admits the parametrization
$$Y(a,\rho, \xi) := \big( \tanh(\rho/2)\xi , \lambda(a,\rho)\big),
~~ \rho \ge a.$$

The geometric data of  $\cC_{a,+}$ are readily calculated. In
particular,

$$|A_a|^2(\rho) = n(n-1) \big( \dfrac{\sinh^{n-1}(a) \cosh(\rho)}
{\sinh^{n}(\rho)} \big)^2,$$ and
$$d\mu_a = \sinh^{2n-2}(\rho) \big( \sinh^{2n-2}(\rho) -
\sinh^{2n-2}(a)\big)^{-1/2} \, d\rho \, d\mu_S.$$

The first assertion follows ($|A_a|^n d\mu_a$ tends to zero
exponentially at infinity). For the second assertion, we use Gauss
equation and minimality to get that

$$K_a = \widehat{K}_a - \frac{1}{2}|A_a|^2 = - v_a^2 - \frac{1}{2}|A_a|^2,$$

where $\widehat{K}_a$ is the sectional curvature of the $2$-plane
tangent to $\cC_a$  in the ambient space $\HH^2\times \R$ and where
$v_a$ is the vertical component of the unit normal to $\cC_a$,
$$v_a(\rho) = \gh(N_a,\partial_t) = \sinh^{1-n}(\rho) \big( \sinh^{2n-2}(\rho)
- \sinh^{2n-2}(a) \big)^{1/2}.$$ Assertion 2 follows because $v_a$
tends to $1$ at infinity on $\cC_{a,+}$. \qed

\subsection{Catenoids in $\HH^n \times \R$, stability properties}
\label{SS-dim3-catsta}

Recall that the catenoid $\cC_a$ is generated by the curve $t
\mapsto \big( \tanh(f(a,t)/2),t \big)$ in the vertical plane $\V$,
where $f(a,\cdot)$ is the maximal solution of the Cauchy problem
(\ref{E-cat3-10}). This yields the parametrization
\begin{equation}\label{E-cat3-12a}
X(a,t,\xi) = \big( \tanh(f(a,t)/2) \xi, t \big)
\end{equation}
for $\cC_a$, with $t \in \R$ and $\xi \in S^{n-1}$. According to
Property \ref{P-frame-1} (2), we have two Jacobi fields on the
catenoid $\cC_a$.

\noib The \emph{vertical Jacobi field} $v(a,t)$ comes from the
vertical translations $(x,t) \mapsto (x,t+\tau)$ in $\HH^n \times
\R$. It is given by $v(a,t) = \gh(N_a, \partial_t)$, where $N_a$ is
the unit normal to $\cC_a$. According to (\ref{E-rot3-7}), it is
given by the formula
\begin{equation}\label{E-cat3-13a}
v(a,t) = f_t(a,t) \big( 1 + f_t^2(a,t) \big)^{-1/2},
\end{equation}
where $f_t$ stands for the derivative with respect to the variable
$t$. Because $t \mapsto f(a,t)$ is even, the function $t \mapsto
v(a,t)$ is odd. \bigskip

\noib The \emph{variation Jacobi field} $e(a,t)$ comes from the
variations with respect to the parameter $a$. It is given by $e(a,t)
= \gh(N_a, \frac{\partial X}{\partial a})$. According to
(\ref{E-rot3-7}) and (\ref{E-cat3-12a}), the function $e(a,t)$ is
given
\begin{equation}\label{E-cat3-14a}
e(a,t) = - f_a(a,t) \big( 1 + f_t^2(a,t) \big)^{-1/2},
\end{equation}
where $f_a$ stands for the derivative with respect to the variable
$a$. Because $t \mapsto f(a,t)$ is even, the function $t \mapsto
e(a,t)$ is even. \bigskip

\noib The Jacobi fields $v(a,t)$ and $e(a,t)$ have nice expressions
when restricted to the upper-half $\cC_{a,+} = \cC_a \cap (\HH^n
\times \R_+)$ of the catenoid $\cC_a$. Indeed, recall that the
function $f(a,\cdot) : [0,T(a)[ \to [0, \infty[$ has an inverse
function $\lambda(a,\rho)$ given by (\ref{E-cat3-10b}). Using the
relationships
$$\lambda\big(a,f(a,t)\big) \equiv t ~~\text{for}~~ t \ge 0$$
and
$$\lambda_{\rho}(a,f)f_t \equiv 1 ~~\text{and}~~ \lambda_a(a,f) +
\lambda_{\rho}(a,f) f_a \equiv 0,$$ we get the following expressions
for $v(a,t)$ and $e(a,t)$ for $t\ge 0$,
\begin{equation}\label{E-cat3-15a}
\left\{%
\begin{array}{lll}
v(a,t) &=& \big( 1 + \lambda_{\rho}^2(a,f(a,t))\big)^{-1/2},\\
&&\\
e(a,t) &=& v(a,t) \lambda_a(a,f(a,t)).\\
\end{array}%
\right.
\end{equation}

For $\rho \ge a$, define the functions $v_1(a,\rho), A_1(a,\rho)$
and $B_1(a,\rho)$, by the following formulas
\begin{equation}\label{E-cat3-16c}
\left\{%
\begin{array}{lll}
v_1(a,\rho) &=& \big( 1 + \lambda_{\rho}^2(a,\rho)
\big)^{-\frac{1}{2}} = \big( \dfrac{\sinh^{2n-2}(\rho) -
\sinh^{2n-2}(a)}{\sinh^{2n-2}(\rho)} \big)^{\frac{1}{2}}\\[6pt]
A_1(a,\rho) &=& \dfrac{\cosh(a)}{\cosh(\rho)} \,
\Big( \dfrac{\sinh(a)}{\sinh(\rho)} \Big)^{n-2},\\[6pt]
B_1(a,\rho) &=& \cosh(a) \int_1^{\frac{\sinh(\rho)}{\sinh(a)}}
(v^{2n-2}-1 )^{-\frac{1}{2}} ( \sinh^2(a) v^2 +1)^{-\frac{3}{2}}\, dv.\\
\end{array}%
\right.
\end{equation}

From (\ref{E-cat3-10b}), we can write
\begin{equation}\label{E-cat3-16a}
\lambda(a,\rho) = \sinh(a) \int_1^{\frac{\sinh(\rho)}{\sinh(a)}}
\big( v^{2n-2}-1 \big)^{-\frac{1}{2}} \big( \sinh^2(a) v^2
+1\big)^{-\frac{1}{2}}\, dv
\end{equation}
and compute $\lambda_a$,
$$
\begin{array}{ll}
\lambda_a(a,\rho) =& - \cosh(a) \sinh^{n-2}(a) \tanh(\rho) \big(
\sinh^{2n-2}(\rho) - \sinh^{2n-2}(a)\big)^{-\frac{1}{2}} +\\[6pt]
& \hphantom{xxx} + \cosh(a) \int_1^{\frac{\sinh(\rho)}{\sinh(a)}} (
v^{2n-2} - 1)^{-\frac{1}{2}}(\sinh^2(a) v^2 + 1)^{-\frac{3}{2}}\, dv.\\
\end{array}
$$

We obtain,
\begin{equation}\label{E-cat3-16d}
\lambda_a(a,\rho) v_1(a,\rho) = - A_1(a,\rho) + B_1(a,\rho)
v_1(a,\rho).
\end{equation}
We summarize the relevant properties in the following lemma whose
proof is straightforward.

\begin{lem}\label{L-cat3-17a}
Define the functions $A(a,t)$ and $B(a,t)$ for $t\ge 0$ by
\begin{equation}\label{E-cat3-17bb}
A(a,t) = A_1(a,f(a,t)), ~~B(a,t) = B_1(a,f(a,t)),
\end{equation}
see Formulas(\ref{E-cat3-16c}). Then,
\begin{equation}\label{E-cat3-17c}
e(a,t) = - A(a,t) + B(a,t) v(a,t), ~~\text{for}~~ t \ge 0.
\end{equation}
Furthermore, for $t \ge 0$,
\begin{enumerate}
    \item $A(a,t) > 0, ~A(a,0)=1$ ~~and~~ $\lim_{t \to T(a)}A(a,t)=0$,
    \item $B(a,t) > 0, ~B(a,0)=0$ ~~and~~ $\lim_{t \to
    T(a)}B(a,t)=C(a)$,
    $\text{where~~} C(a) = \cosh(a) \int_1^{\infty} ( v^{2n-2}-1)^{-\frac{1}{2}}
    ( \sinh^2(a) v^2 +1)^{-\frac{3}{2}}\, dv.$
    \item $v(a,t) = v_1(a,f(a,t))$ for $t>0$, so that\\
    $v(a,t) > 0$ for $t>0$, $v(a,0)=0$ ~and~ $\lim_{t \to
    T(a)}v(a,t)=1$.
\end{enumerate}
\end{lem}

\textbf{Notation}. ~For $\alpha < \beta \in [0,T(a)]$, let
$\cD(\alpha,\beta)$ denote the rotationally symme\-tric domain
\begin{equation}\label{E-cat3-18a}
\cD_a(\alpha,\beta) = X(a,]\alpha,\beta[,S^{n-1}).
\end{equation}
In particular, $\cD_a(0,T(a))$ is the half-vertical catenoid
$\cC_{a,+} = \cC_a \cap (\HH^n \times \R_{+})$.

\begin{thm}\label{T-cat3-20}
The stability properties of the rotationally symmetric domains
$\cD_a(\alpha,\beta)$ on the catenoid $\cC_a$ are as follows.
\begin{enumerate}
    \item There exists some $\sigma(a) \in ]0,T(a)[$ such that the
    relatively compact domain $\cD_a(-\sigma(a),\sigma(a))$ is
    stable-unstable. Hence, for any $\alpha \in ]0,\sigma(a)[$, the domain
    $\cD_a(-\alpha ,\alpha )$ is stable; for any $\alpha \in ]\sigma(a),T(a)[$,
    the domain $\cD_a(-\alpha ,\alpha )$ is unstable.
    \item There exists some $\tau(a) \in ]0,T(a)[$ such that
    \begin{enumerate}
    \item the (non relatively compact) domain $\cD_a(-\tau(a),
    T(a))$ is stable,
    \item for any $\alpha \in ]\tau(a),T(a)[$, there exists some
    $\beta(\alpha) \in ]\tau(a),T(a)[$ such that the domain
    $\cD_a(-\alpha, \beta(\alpha))$ is stable-unstable.
    \end{enumerate}
    \item The catenoid $\cC_a$ has index $1$.
\end{enumerate}
\end{thm}

The above domains are generated by the portions of curves
illustrated in the following figures.

\begin{pb1-figs}
\begin{figure}[h!tb]
\begin{minipage}[l]{.25\linewidth}
    \includegraphics[scale=.22, angle=0]{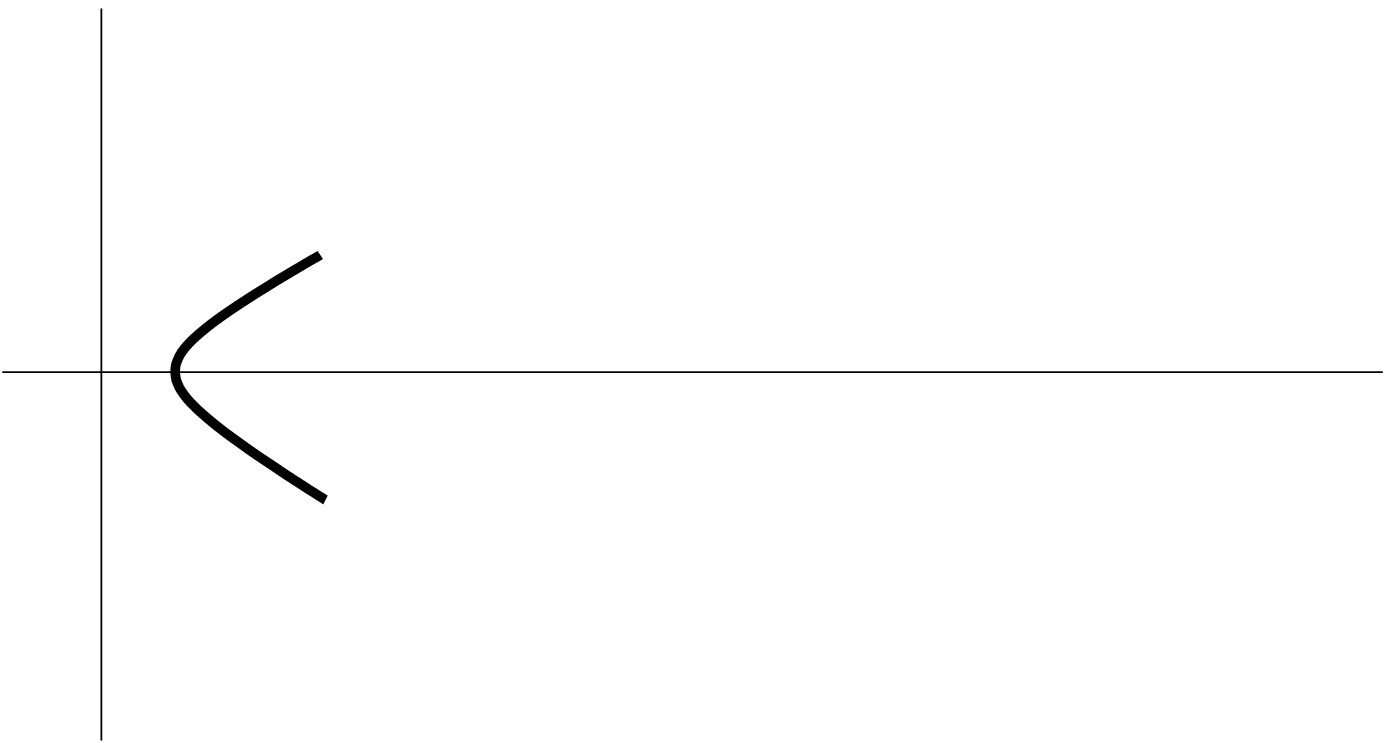}
    \caption[$\cD_a(-\sigma(a),\sigma(a))$]{Case 1}
\end{minipage}\hfill %
\begin{minipage}[r]{.25\linewidth}
    \includegraphics[scale=.22, angle=0]{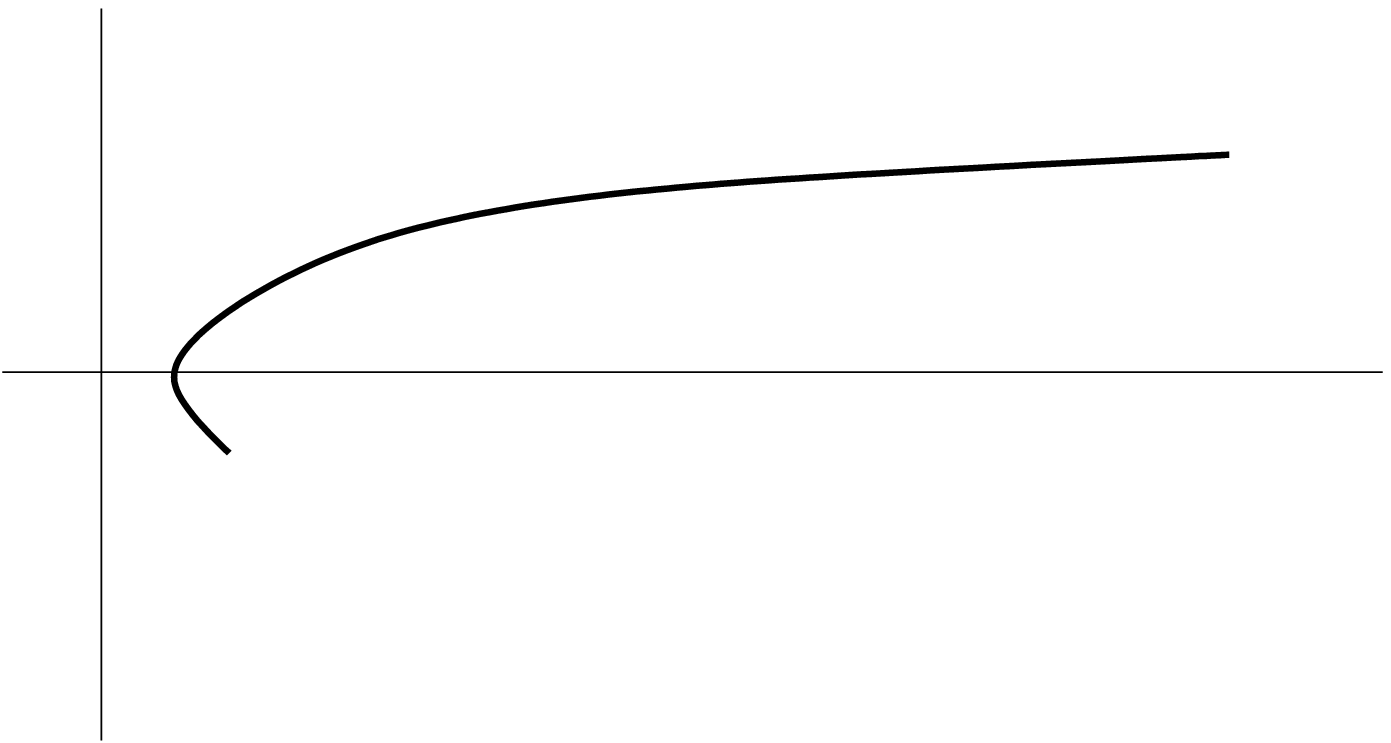}
    \caption[$\cD_a(-\tau(a),T(a))$]{Case 2a}
\end{minipage}\hfill %
\begin{minipage}[l]{.25\linewidth}
    \includegraphics[scale=.22, angle=0]{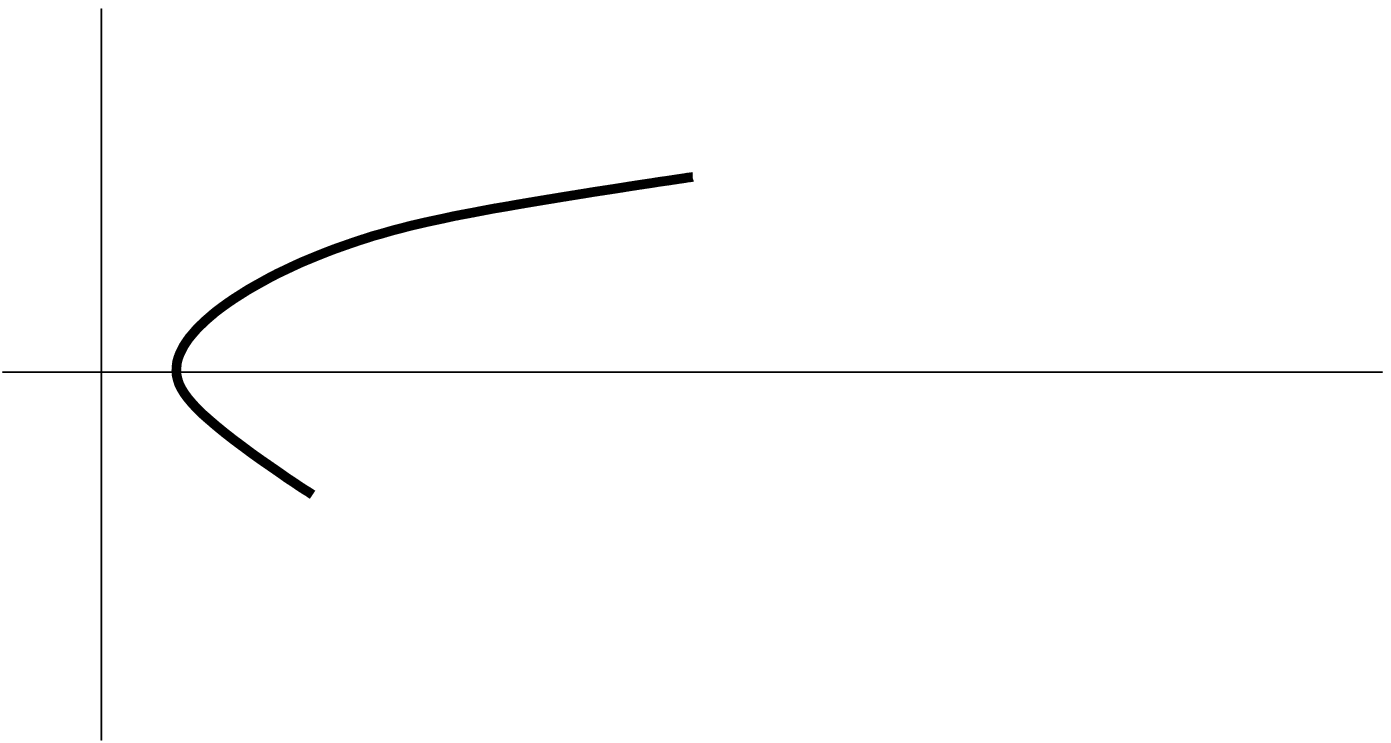}
    \caption[$\cD_a(-\alpha, \beta(\alpha))$]{Case 2b}
\end{minipage}\hfill %
\end{figure}
\end{pb1-figs}

\pf  \emph{Assertion 1}. Consider the function $e(a,t)$. According
to Lemma \ref{L-cat3-17a}, $e(a,0)=-1$ and $\lim_{t \to T(a)}e(a,t)
= C(a) > 0$, so that it must vanish at least once on $]0,T(a)[$. It
turns out (compare with Lemma \ref{L-cat3-21} below) that
$e(a,\cdot)$ has a unique positive zero $\sigma(a)$. Because
$e(a,t)$ is even in $t$, it does not vanish in the open set
$\cD_a(-\sigma(a), \sigma(a))$ and satisfies $J_a(e) = 0
~~\text{in}~~\cD_a(-\sigma(a), \sigma(a))$, and $e|\partial
\cD_a(-\sigma(a), \sigma(a)) = 0.$  This means that
$\cD_a(-\sigma(a), \sigma(a))$ is a stable-unstable domain. The
second assertion follows from Property \ref{P-frame-1} (1). \bigskip

\emph{Assertion 2}. Take any $\alpha \in ]0,T(a)[$ and define the
function $w(a,\alpha,t)$ by
\begin{equation}\label{E-cat3-20a}
w(a,\alpha,t) = e(a,\alpha) v(a,t) + v(a,\alpha) e(a,t),
~~\text{for}~~ t \in ]-T(a),T(a)[.
\end{equation}

This is a Jacobi field on $\cC_a$ and furthermore $w(a,\alpha,
-\alpha) = 0$, because $v$ is odd and $e$ is even with respect to
$t$. Note also that $w(a,\alpha, 0) = - v(a,\alpha) < 0$.

\begin{lem}\label{L-cat3-21}
The function $w(a,\alpha, \cdot)$ vanishes only once on $]-T(a),0[$
and vanishes at most once on $]0,T(a)[$.
\end{lem}

Let us prove the first assertion of the Lemma, the proof of the
second assertion is similar. Assume that $w(a,\alpha, \cdot)$ has at
least two consecutive zeroes $\alpha_1 < \alpha_2$ in the interval
$]-T(a),0[$. The domain $\cD_a(\alpha_1,\alpha_2)$ would then be
stable-unstable because $J_a(w)=0$ on $\cD_a(\alpha_1,\alpha_2)$ and
because $w$ vanishes on $\partial \cD_a(\alpha_1,\alpha_2)$. On the
other-hand, the Jacobi field $v$ satisfies $J_a(v)=0$ and $v < 0$ in
$\cD_a(\alpha_1,\alpha_2)$. By Property \ref{P-frame-1} (3), we have
that $\lambda_1(\cD_a(\alpha_1, \alpha_2)) > 0$ which contradicts
the fact that this domain is stable-unstable. This proves the lemma.
\bigskip

In order to determine whether the function $w(a,\alpha, \cdot)$
vanishes on $]0,T(a)[$ or not, it is sufficient to look at the
behaviour of $w(a,\alpha, t)$ when $t$ tends to $T(a)$ from below.
For this purpose, we use the expression (\ref{E-cat3-17c}) for
$e(a,t)$ and we write
$$w(a,\alpha, t) = - A(a,t) v(a,\alpha) + v(a,t) \big( e(a,\alpha) +
B(a,t) v(a,\alpha)\big).$$

Using Lemma \ref{L-cat3-17a}, we can write
$$W(a,\alpha) := \lim_{t \to T(a)}w(a,\alpha, t) = e(a,\alpha) +
C(a) v(a,\alpha).$$

If $W(a,\alpha) \le 0$, then $w(a,\alpha, t)$ does not vanish on
$]0,T(a)[$ and in fact on $]-\alpha,T(a)[$; if $W(a,\alpha) > 0$,
then $w(a,\alpha, t)$ has one and only one zero $\beta(\alpha)$on
$]0,T(a)[$. \bigskip

We now observe that $W(a,t) := e(a,t) + C(a) v(a,t)$ is a Jacobi
field on $]0,T(a)[$ which take the value $-1$ at $0$ and the value
$C(a) v(a,\sigma(a)) > 0$ at $\sigma(a)$. It follows from Lemma
\ref{L-cat3-21} that $W(a,\cdot)$ has one and only one positive zero
$\tau(a) \in ]0,\sigma(a)[$. We have that $W(a,t) \le 0$ on
$]0,\tau(a)[$, so that for any $\alpha \in ]0,\tau(a)]$, the
function $w(a,\alpha, t)$ has only one zero $-\alpha$ on
$]-T(a),T(a)[$. This proves the Assertion 2(a). On the other-hand,
$W(a,t) > 0$ on $]\tau(a),T(a)[$, so that for any $\alpha \in
]\tau(a),T(a)[$, the function $w(a,\alpha, t)$ has a unique positive
zero $\beta(\alpha) \in ]0,T(a)[$. This proves the Assertion 2(b).
\bigskip

\emph{Assertion 3}. Assertion 1 shows that $\cC_a$ has index at
least $1$. In order to show that the index is at most one, we use
Fourier decomposition with respect to the variable $\xi$ and an
extra stability argument. \bigskip

Recall that we work in the ball model for $\HH^n$. Let $\gamma$ be a
geodesic through $0$ in $\HH^n$. Up to a rotation, we may assume
that $\gamma(s) = \big( \tanh(s/2), 0, \cdots, 0 \big)$. Let
$\HH^n_{+} = \{(x_1,\cdots,x_n) \in \B ~|~ x_1>0\}$ and let
$\cC_{a,\gamma +} = \cC_a \cap (\HH^n_+ \times \R)$. We call this
set a \emph{half-horizontal catenoid}. \bigskip

\emph{Claim 1}. A half-horizontal catenoid $\cC_{a,\gamma +}$ is
stable. \bigskip

To prove the claim, we shall find a positive Jacobi field on
$\cC_{a,\gamma +}$.

Let $z=x+iy$ denote the complex coordinate in $\HH^2$ (ball model).
We consider the group of hyperbolic isometries along the geodesic
$\gamma$ and we extend these isometries slice-wise as isometries in
$\HH^2\times \R$. We then have the one-parameter group of isometries
$$(z;t) \mapsto \big( \dfrac{e^{\tau}(1+z) - (1-z)}{e^{\tau}(1+z)
+ (1-z)} ; t\big) \text{~~in~~} \HH^2\times \R.$$

The associated Killing vector-field in $\HH^2 \times \R$ is given by
$\cK_{\gamma}(z;t) = \big(\frac{1}{2}(1-z^2) ; 0 \big)$ or, in the
$(x,y)$ coordinates, $\cK_{\gamma}(x,y;t) =
\big(\frac{1}{2}(1-x^2+y^2), -xy ; 0 \big)$ which can be written as

$$\cK_{\gamma}(x,y;t) = \frac{1}{2}\big(1+x^2+y^2)\big)(1,0;0) - x(x,y;0)$$

where $(1,0;0)$ and $(x,y;0)$ are seen as vectors in $\R^2 \times \R
= T_{(x,y;t)}\HH^2\times \R$. \bigskip

This formula can easily be generalized to higher dimensions as

$$\cK_{\gamma}(x;t) = \frac{1}{2}(1+|x|^2 (e_1;0) - x_1(x;0),$$

where $x = (x_1, \cdots, x_n), e_1=(1, 0,\cdots,0), |x|^2 = x_1^2 +
\cdots + x_n^2$, and where $(e_1;0)$ and $(x;0)$ are seen as vectors
in $\R^n \times \R = T_{(x,t)}\HH^n \times \R$. Writing the point
$x$ in the parametrization $X$ as $x = \tanh(f(a,t)/2)\xi$, we
obtain that

$$\cK_{\gamma}(\tanh(f(a,t)/2)\xi;t) = \frac{1}{2} \big( 1 + \tanh^2(f/2)\big)
(e_1;0) - \tanh^2(f/2) \xi_1 (\xi;0).$$

Using the fact that $(1+f_t^2)^{-1/2} = \big(
\frac{\sinh(a)}{\sinh(f)} \big)^{n-1}$ on $\cC_a$, we find that the
Killing field $\cK_{\gamma}$ gives rise to the horizontal Jacobi
field

$$h_{\gamma}(a,t,\xi) = \big( \dfrac{\sinh(a)}{\sinh(f(a,t))} \big)^{n-1}\xi_1$$

which is positive on $\cC_{a,\gamma +}$. \bigskip

\emph{Claim 2}. On $S^{n-1}$ equiped with the standard Riemannian
metric, there exists an orthonormal basis of spherical harmonics
$Y_k, k\ge 0$ with the property that the nodal domains of the $Y_k,
k\ge 1$ are contained in hemispheres. \bigskip

The property is clearly true on $S^1$ and can be proved by induction
on the dimension, using polar coordinates centered at a given point
on the sphere. \bigskip

\emph{Claim 3}. The Jacobi operator on $\cC_a$ can be written as
$$J_a = L_{a,t} - q(a,t) \Delta_{S,\xi},$$
where $L_{a,t}$ is a Sturm-Liouville operator on the $t$ variable,
with coefficients depending only on $a$ and $t$, where $q(a,t)$ is a
positive function and where $\Delta_{S,\xi}$ is the Laplacian of the
sphere $S^{n-1}$ acting on the $\xi$-variable. This claim follows
immediately from Formulas (\ref{E-rot3-5}) and (\ref{E-rot3-6}) for
the metric and the Riemannian measure on a rotation hypersurface and
from the expression for the quadratic form associated with $J_a$.
\bigskip

Assume that the index of $\cC_a$ is at least $2$. Then, there exists
some $S \in ]0,T(a)[$ such that $J_a$ has at least two negative
eigenvalues $\lambda_1(S) < \lambda_2(S) < 0$ in $\cC_a(-S,S)$ (we
only consider Dirichlet boundary conditions). Because the least
eigenvalue $\lambda_1(S)$ is simple, a corresponding eigenfunction
$u$ must be invariant under rotations (\ie only depends on the
variable $t$) and say positive. Consider an eigenfunction $v$
associated with $\lambda_2(S)$. We claim that $v$ cannot be
invariant under rotations. Indeed, $v$ would otherwise depend only
on the variable $t$, it would be orthogonal to $u$ and hence it
would have to vanish on $]-S,S[$. This would contradict the fact
that the domains $\cC_a(-S,0)$ and $\cC_a(0,S)$ are stable. Since
$v$ is not rotationally invariant, there exists some $p \ge 1$ and
some $v_p \not = 0$ in the decomposition into spherical harmonics
with respect to the second variable, $v(t,\xi) = \sum_{k=0}^{\infty}
v_k(t) Y_k(\xi)$. We would have $J_a(v_p Y_p) = \lambda_2(S) v_p
Y_p$. Using Claim 2 and the fact that $\lambda_2(S) < 0$, this would
mean that any nodal domain of $v_p Y_p$ is unstable, in
contradiction with Claim 1.
\bigskip

Assuming that the index is at least $2$ therefore yields a
contradiction and hence the index of $\cC_a$ is exactly one. \qed
\bigskip

\textbf{Remark}. ~ It follows from the positivity of the Jacobi
field $v(a,t)$ for $t \in ]0,T(a)[$ that the upper half-catenoid
$\cC_{a,+}$ is stable (in the sense that any relatively compact
domain $\Omega$ contained in $\cC_{a,+}$ is stable, see Section
\ref{SS-jac}). The second assertion in the preceding theorem says
more. Indeed, there exists some $\tau(a) \in ]0,T(a)[$ such that the
non-compact domain $\cD_a(-\tau(a), T(a))$ is stable and stricly
contains $\cC_{a,+}$. This is different from what happens for
Euclidean catenoids. Indeed, the half-catenoid $\cC_{a,+}$ in $\R^3$
is a \emph{maximal} stable domain (\cite{Lin870}). We study this
phenomenon with more details in \cite{BSA09L}. \bigskip

\begin{pb1-figs}
\begin{figure}[h!tb]
\begin{minipage}[l]{.4\linewidth}
    \includegraphics[scale=.25, angle=-90]{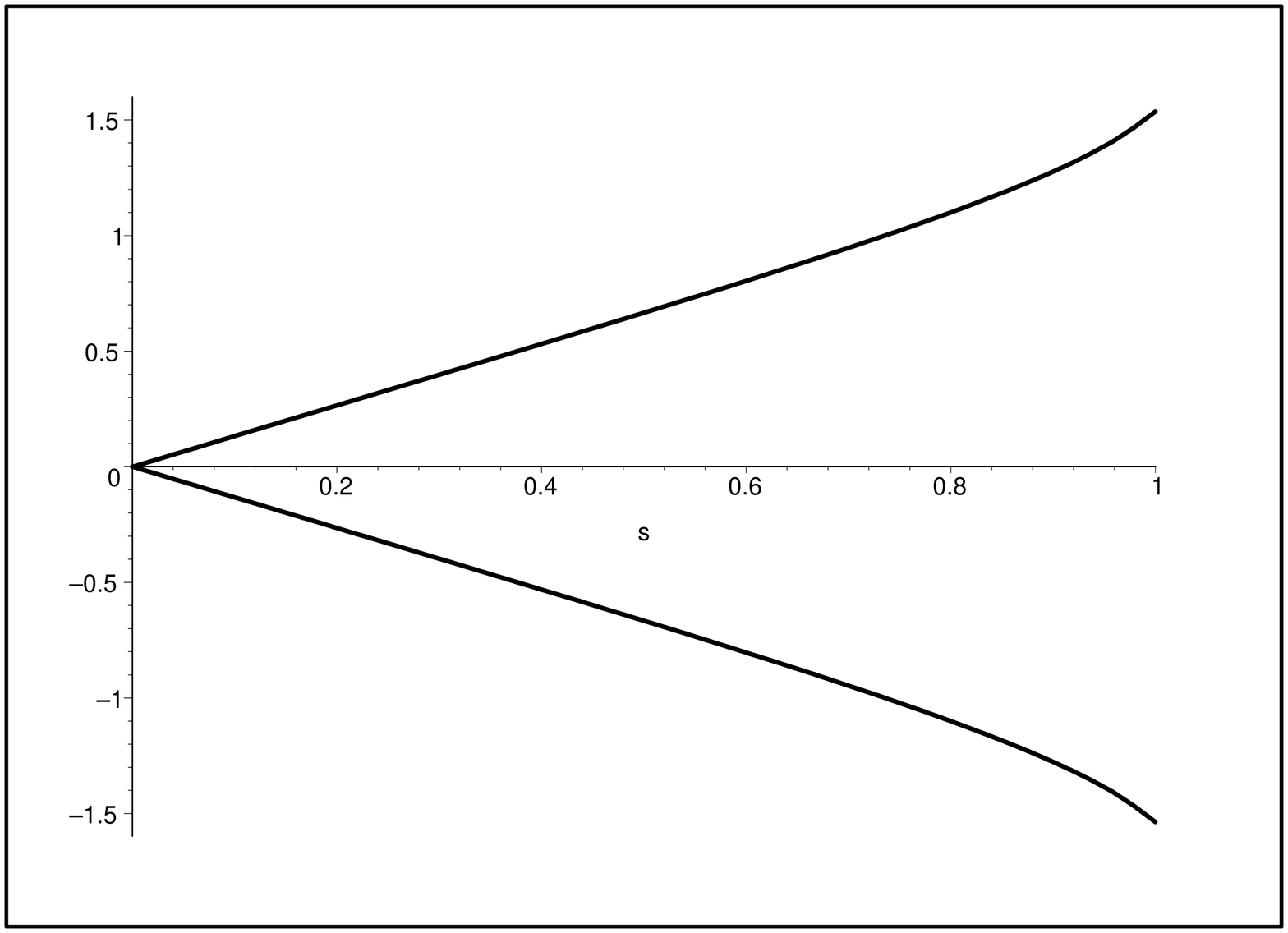}
    \caption[Envelope, $n=2$]{Envelope, $n=2$}
    \label{F-envelope}
\end{minipage}\hfill %
\begin{minipage}[r]{.4\linewidth}
    \includegraphics[scale=.25, angle=-90]{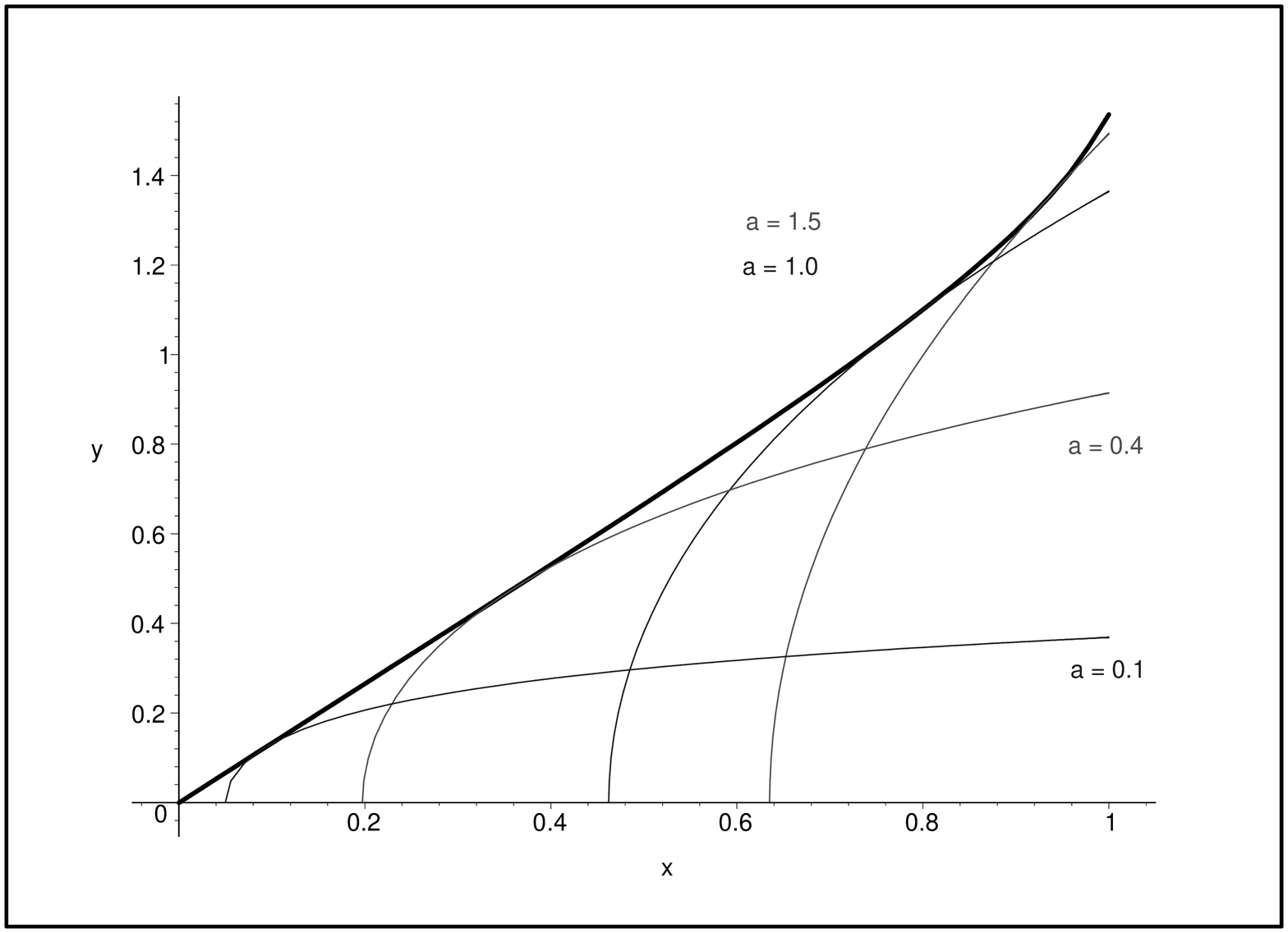}
    \caption[Catenaries and envelope, $n=2$]{Catenaries, envelope}
    \label{F-catene}
   \end{minipage}
\end{figure}
\end{pb1-figs}

\textbf{Geometric interpretation}. ~According to Proposition
\ref{P-cat3-10a}, Assertion (6), two distinct catenaries $C_a$ and
$C_b$ meet at exactly two symmetric points, $m_{\pm}(a,b)$. Fixing $a$
and letting $b$ tend to $a$, the points $m_{\pm}(a,b)$ tend to limit
points $m_{\pm}(a)$ which correspond to the points where the catenary
$C_a$ touches the envelope of the family of catenaries
$\{C_a\}_{a>0}$. According to \cite{V50}, \S 58, page 127 ff, the
condition defining the envelope a family $\Gamma_a$, given by the
parametrization $\big( x(a,t), y(a,t) \big)$, is the condition
$$\begin{vmatrix}
  x_a(a,t) & x_t(a,t) \\
  y_a(a,t) & y_t(a,t) \\
\end{vmatrix} = 0.$$

Specializing to catenaries, we find that the envelope condition is
precisely the condition that $e(a,t)=0$. Therefore, the value
$\sigma(a)$ is precisely the value of $t$ at which the catenary
$C_a$ touches the envelope of the family, see Figure~\ref{F-catene}.

\begin{cor}\label{C-cat3-21}
The stable-unstable do\-main $\cD_a(-\sigma(a),\sigma(a))$ is
preci\-sely the sym\-metric, rotation invariant compact domain
bounded by the two spheres where the catenoid $\cC_a$ touches the
envelope of the family.
\end{cor}

\subsection{Translation invariant hypersurfaces in $\HH^n \times \R$}\label{SS-dim3-trans1}

Let $\gamma$ be a complete geodesic through $0$ in the ball model
$\B$ of the hyperbolic space $\HH^n$, parametrized by the signed
distance $\rho$ to $0$. Let $\PP$ be the hyperbolic hyperplane
orthogonal to $\gamma$ at $0$. We consider the hyperbolic
translations along the geodesics passing through $0$ in $\PP$. The
image of a point of $\gamma$ under these translations is an
equidistant hypersurface to $\PP$ in $\HH^n$. We can extend these
translations ``slice-wise" to give positive isometries of $\HH^n
\times \R$ which we call \emph{hyperbolic translations}. \bigskip

A generating curve $(\tanh(\rho /2), \mu(\rho))$ in the vertical
Euclidean plane $\gamma \times \R$ gives rise, under the previous
isometries, to a translation invariant hypersurface $M
\hookrightarrow \HH^n \times \R$, whose intersection with the slice
$\HH^n \times \{\mu(\rho)\}$ is the equidistant hypersurface to $\PP
\times \{\mu(\rho)\}$ in the slice, at distance $\rho$. \bigskip

The principal directions of curvature of $M$ are the tangent vector
to the generating curve and the directions tangent to the
equidistant hypersurface. The corresponding principal curvatures are
given respectively by
\begin{equation}\label{E-trans3-1E}
k_G(\rho) = \ddot{\mu}(\rho) \big( 1 +
\dot{\mu}^2(\rho)\big)^{-3/2},
\end{equation}
and
\begin{equation}\label{E-trans3-1G}
k_E(\rho) = \dot{\mu}(\rho) \big( 1 + \dot{\mu}^2(\rho)\big)^{-1/2}
\, \tanh(\rho).
\end{equation}

\begin{pb2}
\como{PB2}

\input{pb2-03}

\comf{PB2}\bigskip
\end{pb2}

It follows that the mean curvature of $M$ is given by
\begin{equation*}\label{E-trans3-2}
n H(\rho)  = \ddot{\mu}(\rho) \big( 1 +
\dot{\mu}^2(\rho)\big)^{-3/2} + (n-1) \dot{\mu}(\rho) \big( 1 +
\dot{\mu}^2(\rho)\big)^{-1/2} \, \tanh(\rho)
\end{equation*}
or, equivalently, by
\begin{equation}\label{E-trans3-3}
n H(\rho) \cosh^{n-1}(\rho) = \partial_{\rho} \Big( \dot{\mu}(\rho)
\big( 1 + \dot{\mu}^2(\rho)\big)^{-1/2} \, \cosh^{n-1}(\rho)\Big).
\end{equation}

This formula allows us to study constant mean curvature
hypersurfaces invariant by hyperbolic translations. In this paper we
only consider the case $H=0$ and we refer to \cite{BS08b} for the
case $H \not = 0$.

\subsection{Translation invariant minimal hypersurfaces in
$\HH^n \times \R$}\label{SS-dim3-trans2}

In this section, we establish the following theorem which
generalizes the $2$-dimen\-sio\-nal result of \cite{ST08}.

\begin{thm}\label{T-trans3-1}
There exists a $1$-parameter family $\{ \cM_d,\, d>0\}$ of complete
embedded minimal hypersurfaces in $\HH^n \times \R$ invariant under
hyperbolic translations. The hypersurfaces are $\cM_d$ stable (in
the sense of the Jacobi operator), their principal curvatures go
uniformly to
zero at infinity, but they have infinite total curvatures. 

More precisely,
\begin{enumerate}
    \item If $d>1$, the hypersurface $\cM_d$ consists of the union of
    two symmetric vertical graphs over the exterior of an equidistant
    hypersurface in the slice $\HH^n\times \{0\}$. It is also a horizontal
    graph, and hence stable.\\[2pt]
    The family $\cM_d$ has finite vertical height $h_T(d)$, a function which
    decreases from infinity to $\pi /(n-1)$. In particular, it is bounded
    from below by $\pi/(n-1)$, the upper bound of the heights of the family of
    catenoids. Furthermore, the asymptotic boundary of $\cM_d$  consists of the
    union of two copies of an hemisphere $S_+^{n-1}\times \{0\}$ of $\partial_\infty
    \HH^n \times \{0\}$ in parallel slices $t=\pm S(d)$, glued with the finite
    cylinder $\partial S_+^{n-1}\times [-S(d),S(d)]$.
    \item If $d=1$, the hypersurface $\cM_1$ is a complete stable vertical
    graph over a half-space in $\HH^n\times \{0\}$, bounded by a totally
    geodesic hyperplane $P$. It takes infinite boundary value data on
    $P$ and constant asymptotic boundary value data. Furthermore, the asymptotic
    boundary of $\cM_1$ is the union of a spherical
    cap $S$ in $\partial_\infty \HH^n\times \{c\}$ with a half-vertical cylinder
    over $\partial S$.
    \item If $d<1$, the hypersurface $\cM_d$ is an entire stable vertical graph
    with finite vertical height.
    Furthermore, its asymptotic boundary consists of a homologically
    non-trivial $(n-1)$-sphere in $\partial_\infty \HH^n \times \R.$
\end{enumerate}
\end{thm}

\pf  In the minimal case, Equation (\ref{E-trans3-3}) can be written
\begin{equation}\label{E-trans3-4}
\dot{\mu}(\rho) \big( 1 + \dot{\mu}^2(\rho)\big)^{-1/2} \,
\cosh^{n-1}(\rho) = d,
\end{equation}
for some constant $d$ which satisfies $d \le \cosh^{n-1}(\rho)$ for
all $\rho$ for which the solution exists. Changing $\mu$ to $- \mu$
if necessary, we may assume that $d$ is non-negative and hence,
$\dot{\mu}(\rho) \ge 0$ and $\dot{\mu}(\rho) = d \big(
\cosh^{2n-2}(\rho) - d^2\big)^{-1/2}$ whenever the square root
exists. We have to consider three cases, $d > 1$, $d=1$ and $d<1$.

\bigskip
$\boxed{d>1}$ ~~ Let $d =: \cosh^{n-1}(a)$, with $a>0$. It follows
from Equation (\ref{E-trans3-4}) that
$$\dot{\mu}(\rho) = \cosh^{n-1}(a) \big( \cosh^{2n-2}(\rho) - \cosh^{2n-2}(a)
\big)^{-1/2}.$$
Up to a vertical translation, the solution $\mu_{+}(a,\rho)$ of
Equation (\ref{E-trans3-4}) is given by
\begin{equation}\label{E-trans3-5}
\mu_{+}(a,\rho) = \cosh^{n-1}(a) \int_a^{\rho} \big( \cosh^{2n-2}(r)
- \cosh^{2n-2}(a) \big)^{-1/2} \, dr
\end{equation}
or, making $\cosh(r) = \cosh(a) \, t$,
\begin{equation}\label{E-trans3-6}
\mu_{+}(a,\rho) = \cosh(a) \, \int_1^{\frac{\cosh(\rho)}{\cosh(a)}}
(t^{2n-2} - 1)^{-\frac{1}{2}} \, (\cosh^2(a) t^2 - 1)^{-\frac{1}{2}}
\, dt.
\end{equation}

These integrals converge at $\rho = a$ (\resp at $t =1$) and at
infinity and
\begin{equation}\label{E-trans3-7}
h_T(d) := 2 \cosh(a) \, \int_1^{\infty} (t^{2n-2} - 1)^{-1/2} \,
(\cosh^2(a) t^2 - 1)^{-1/2} \, dt
\end{equation}

is the height of the hypersurface $\cM_{\cosh^{n-1}(a)}$. The
function $h_T(d)$ is decreasing in $d$, tends to infinity when $d$
tends to $1_+$ and to $\pi/(n-1)$ when $d$ tends to infinity.
(Hints. When $a$ tends to zero, use the fact that (\ref{E-trans3-7})
is bigger than some constant times the integral $\int_1^2 \big(
(t-1) (\cosh(a) t - 1) \big)^{-1/2} \, dt$ which can be computed
explicitly. When $a$ tends to infinity, use the fact that $\int (t^N
- 1)^{-1/2}t^{-1}\, dt = \frac{2}{N} \arctan \sqrt{t^N-1}$.) The
assertions on the asymptotic boundary are clear.

\begin{figure}[h!tb]
\begin{minipage}[l]{.42\linewidth}
    \includegraphics[scale=0.27, angle=-90]{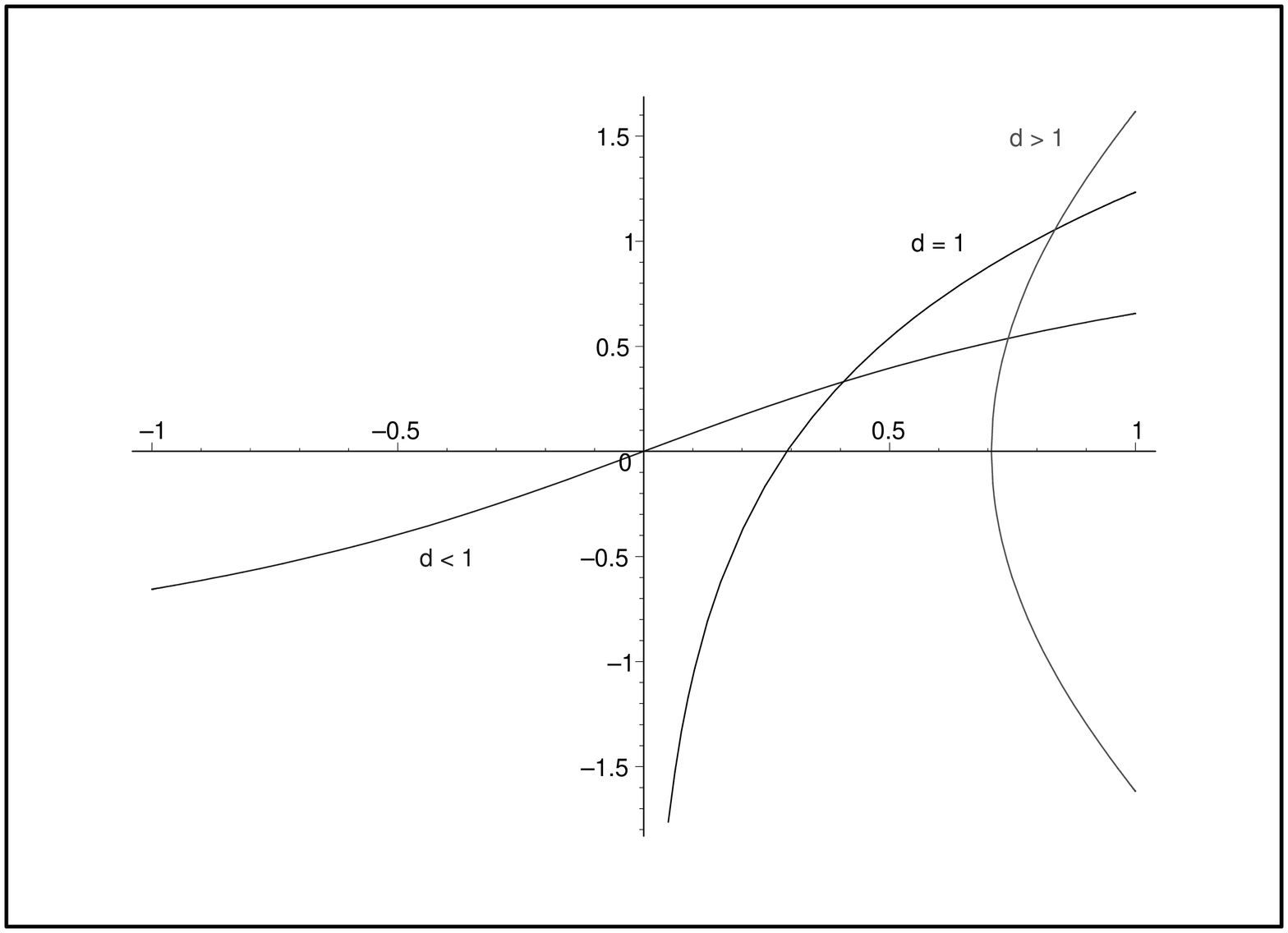}
    \caption[Generatrices for translation invariant hypersurfaces, $n=2$]
    {Generatrices translation invariant hypersurfaces, $n=2$}
    \label{F-dim3-tra2}
\end{minipage} \hfill
\begin{minipage}[r]{.42\linewidth}
    \includegraphics[scale=0.27, angle=-90]{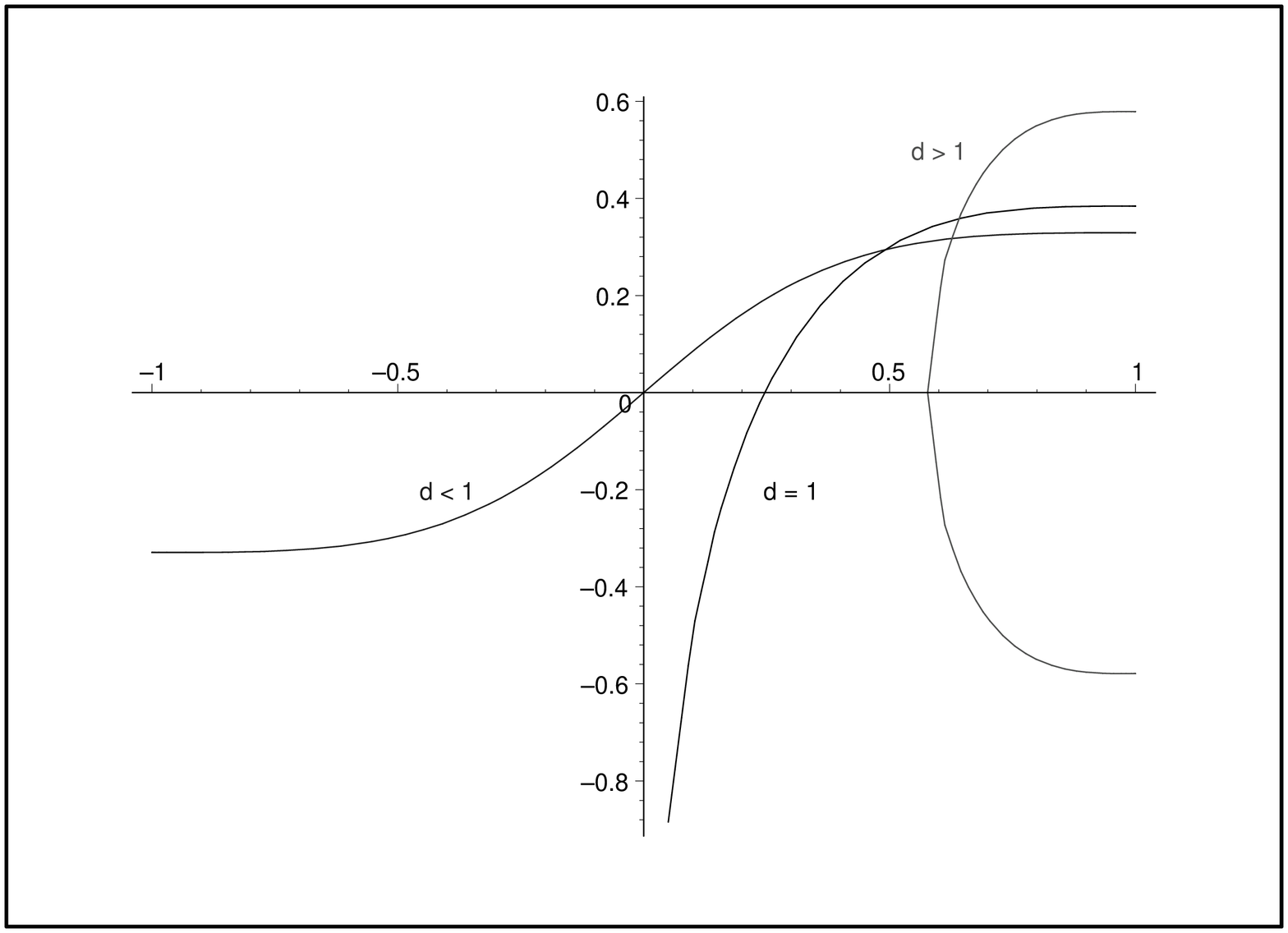}
    \caption[Generatrices for translation invariant hypersurfaces, $n=4$]
    {Generatrices translation invariant hypersurfaces, $n=4$}
    \label{F-dim3-tra4}
\end{minipage}
\end{figure}

$\boxed{d=1}$ ~~ It follows from Equation (\ref{E-trans3-4}) that
$$\dot{\mu}(\rho) = \big( \cosh^{2n-2}(\rho) - 1\big)^{-1/2},$$
so that, when $d=1$, the solution is given by
\begin{equation}\label{E-trans3-10}
\mu_{0}(\rho) = \int_b^{\rho} \big( \cosh^{2n-2}(r) - 1\big)^{-1/2}
\, dr,
\end{equation}
for some constant $b > 0$, and $\mu_0(\rho)$ tends to $- \infty$
when $\rho$ tends to zero and to a finite value when $\rho$ tends to
infinity. The corresponding hypersurface is complete. It is a
vertical graph so that it is stable. The assertion on the asymptotic
boundary is clear.

\bigskip
$\boxed{0< d<1}$ ~~ In this case, Equation (\ref{E-trans3-4}) gives
the following solution (up to a vertical translation),
\begin{equation}\label{E-trans3-15}
\mu_{-}(d,\rho) = d \int_0^{\rho} \big( \cosh^{2n-2}(r) -
d^2\big)^{-1/2} \, dr.
\end{equation}

The corresponding curve can be extended by symmetry and we get a
complete hypersurface in a vertical slab with finite height. This
surface is an entire vertical graph (hence stable). The assertion on
the asymptotic boundary is clear. \qed \bigskip

The generating curves for translation invariant minimal
hypersurfaces are given in Fig. \ref{F-dim3-tra2} and
\ref{F-dim3-tra4}. Note that they cannot meet tangentially at finite
distance.\bigskip

\textbf{Remark.}~ Using the catenoids and the minimal translations
hypersurfaces, the second author and E. Toubiana have extended the
$2$-dimensional results of their paper \cite{ST08} to higher
dimensions, see \cite{ST09}.

\section{Index and total curvature for minimal hypersurfaces in
$\HH^n \times \R$}\label{S-index}

\subsection{Dimension two, $M^2 \looparrowright \HH^2 \times \R$}

For oriented minimal surfaces in $\HH^2 \times \R$ we have the
following general theorem in which we consider two possible notions
of total curvature.

\begin{thm}\label{T-dim2}
Let $M \looparrowright \HH^2 \times \R$ be a complete oriented
minimal immersion with unit normal field $N_M$. Let $v_M :=
\hat{g}(N_M, \partial_t)$ be the vertical component of $N_M$, let
$A_M$ be the second fundamental form of $M$ and let $K_M$ be the
intrinsic curvature of $M$.
\begin{enumerate}
    \item If the integral $\int_M |A_M|^2 \, d\mu_M$ is finite, then
    $A_M$ tends to zero uniformly at infinity.
    \item If the integral $\int_M |K_M| \, d\mu_M$ is finite, then
    $A_M, v_M$ and $K_M$ tend to zero uniformly at infinity.
    \item If the integral $\int_M |K_M| \, d\mu_M$ is finite, then
    the Jacobi operator of $M$ has finite index.
\end{enumerate}
\end{thm}

\textbf{Remarks.}
\begin{enumerate}
    \item For complete orientable minimal surfaces in $\R^3$, finiteness of the index
    is equi\-va\-lent to finiteness of the intrinsic curvature (see
    \cite{FC82, CS90, LR98}).
    No such statement can hold in $\HH^2 \times \R$. Indeed, the surfaces $\cM_d$
    (\cite{ST08} and Section~\ref{SS-dim3-trans2}) are stable complete minimal
    surfaces, invariant under a group of hyperbolic translations. Their total
    curvatures are infinite, so that the converse to Assertion 3 is false.
    \item The assumption $\int_M |K_M| \, d\mu_M$ finite is natural in view of
    Huber's theorem. In \cite{HR08}, L. Hauswirth and H. Rosenberg
    show that this assumption actually implies that the total
    intrinsic curvature is an integer multiple of $2\pi$. There are actually many
    examples of such surfaces (\cite{CR07, HR08}).
    \item As pointed out in the introduction, Assertion 2 is contained in
    \cite{HR08}, Theorem 3.1 whose proof actually gives a $C^2$-control
    on the curvature at infinity. We provide a simple proof of Assertion~2 for
    completeness.
\end{enumerate}

One can slightly improve the above theorem with the following
proposition.

\begin{prop}\label{P-dim2}
The notations are the same as in Theorem \ref{T-dim2}.
\begin{enumerate}
    \item Assume that $5 v_M^2 \le 1$. Then there exists a universal
    constant $C$ such that if the integral $\int_M |A_M|^2 \, d\mu_M$
    is less than $C$ then $M$ is a vertical plane.
    \item Assume that the integral $\int_M |A_M|^2 \, d\mu_M$ is finite
    and that there exists a compact set $\Omega \subset M$
    and a positive constant $c$ such that $v_M^2 \le 1-c < 1$ on
    $M \setminus \Omega$. Then the Jacobi operator of the immersion has finite index.
\end{enumerate}
\end{prop}

\textbf{Remarks.}
\begin{enumerate}
    \item Assertion 1 generalizes the following facts : (i) A
    minimal surface whose intrinsic curvature $K_M$ is zero is part of a
    vertical plane $\gamma \times \R$ (where $\gamma$ is a geodesic
    in $\HH^2$). Indeed, we have $K_M = - \frac{1}{2}|A_M|^2 -
    v_M^2$ and hence $M$ is totally geodesic with horizontal normal
    vector. (ii) A complete minimal surface whose total intrinsic
    curvature is less than $2 \pi$ is a vertical plane (see \cite{HR08}).
    \item We do not know whether the sole assumption $\int_M |A_M|^2 \, d\mu_M$
    finite is sufficient to insure the finiteness of the index of
    the Jacobi operator of $M$.
\end{enumerate}

\textbf{Proof of Theorem \ref{T-dim2} and Proposition
\ref{P-dim2}.}\medskip

\textbf{Fact 1.}~ The function $u := |A_M|$ satisfies the non-linear
elliptic inequality
\begin{equation}\label{E-dim2-sim3}
    - u \, \Delta_M u \le u^4 + (5 \Kh_M +1) u^2 \le u^4 + 4 u^2,
\end{equation}

where $\Kh_M$ is the sectional curvature of the $2$-plane $TM$ in
$\HH^2 \times \R$.

\begin{pb3}
\como{PB3}

\input{pb3-03}

\comf{PB3}\bigskip
\end{pb3}

This formula follows from J. Simons' equation for minimal
submanifolds (\cite{Si68}), from our context ($\HH^2 \times \R$ is
locally symmetric and we work in codimension $1$), and from explicit
curvature computations in $\HH^2 \times \R$.
\bigskip

\textbf{Fact 2.}~ The surface $M$ satisfies the Sobolev inequality
\begin{equation}\label{E-dim2-sob1}
\| f \|_{2} \le S(M) \| df \|_1
\end{equation}

for some positive constant $S(M)$ and all $C^1$ functions $f$ with
compact support. This follows from \cite{HS74, HS75} and the fact
that $\HH^2 \times \R$ is simply-connected and non-positively
curved.\bigskip

\textbf{Fact 3.}~ Curvature computations in $\HH^2 \times \R$ give
\begin{equation}\label{E-dim2-curv}
\Kh_M = - v_M^2 \text{ ~and~ } \rich(N_M,N_M) = - (1-v_M^2).
\end{equation}

The fact that $v_M$ is a Jacobi field implies that
\begin{equation}\label{E-dim2-jac}
    - \Delta_M v_M = v_M^3 + (|A_M|^2 -1) v_M
\end{equation}

and a similar inequality for $|v_M|$.\bigskip

\textbf{Theorem, Assertion 1.}~ Following the general ideas of
\cite{SSY75}, we use (\ref{E-dim2-sim3}) and (\ref{E-dim2-sob1}), to
estimate the $L^p$-norms of $u$ and the classical de
Giorgi-Moser-Nash method to estimate $\| u \|_{\infty}$ outside big
balls. The details appear in the proof of Theorem 4.1, p. 282 of
\cite{BCS98}, where it is observed that the proof only uses the
facts that $u$ satisfies Simons' inequality and $M$ a Sobolev
inequality.

\textbf{Theorem, Assertion 2.}~ By Gauss equation, the Gauss
curvature $K_M$ of $M$ satisfies
\begin{equation}\label{E-dim2-gauss}
    K_M = - \frac{1}{2}|A_M|^2 - v_M^2.
\end{equation}

The assumption implies that both integrals $\int_M |A_M|^2 \,
d\mu_M$ and $\int_M v_M^2 \, d\mu_M$ are finite. By Assertion 1, we
already know that $|A_M|$ tends to zero at infinity, and hence that
it is bounded. Equation (\ref{E-dim2-jac}) then tells us that
$|v_M|$ satisfies an elliptic inequality similar to
(\ref{E-dim2-sim3}) and we can again apply the de Giorgi-Moser-Nash
method to conclude.\bigskip

\textbf{Theorem, Assertion 3 and Proposition, Assertion 2.}~
According to Section \ref{SS-jac} and to the above curvature
calculations, the Jacobi operator can be written as $J_M = -\Delta_M
+ 1 -|A_M|^2 - v_M^2$.  We now follow \cite{BCS97}, Section $2$. It
follows from Assertion 2 in the Theorem that $J_M$ is bounded from
below, essentially self-adjoint and that its essential spectrum lies
above $1$. Because the eigenvalues below the essential spectrum can
only accumulate at $- \infty$ or at the bottom of the essential
spectrum, it follows that $J_M$ has finite index.\bigskip

\textbf{Proposition, Assertion 1.}~ It is a classical fact
(\cite{CP80} and \cite{SSY75}) that Simons' inequality
(\ref{E-dim2-sim3}) can be improved to
$$|du|^2-u \, \Delta_M u \le u^4 + (5 \Kh_M +1) u^2.$$

We then use the expression of $\Kh_M$ and the assumption on $v_M$ to
obtain
$$(a)~~~ |du|^2-u \, \Delta_M u \le u^4.$$

Multiply equation (a) by some function with compact support $\xi$
(to be chosen later on) and integrate by parts to obtain,
$$2 \int_M \xi^2 |du|^2 + 2 \int_M \xi u \bra du,d\xi \ket \le \int_M \xi^2 u^4.$$
Using Cauchy-Schwarz inequality, we obtain
$$(b)~~~ \int_M \xi^2 |du|^2 \le \int_M \xi^2 u^4 + \int_M u^2
|d\xi|^2.$$

Plug the function $f = \xi u^2$ into Sobolev inequality
(\ref{E-dim2-sob1}) to obtain
$$\int_M \xi^2 u^4 \le S \big(\int_M |d(\xi u^2)| \big) \le
2S \big( \int_M u^2 |d\xi|\big)^2 + 8S (\int_M
u^2) \big( \int_M u \xi  |du|\big)^2,$$ where we have noted $S$ for
$S(M)$ and used the fact that $\int_M u^2$ is finite. Using
Cauchy-Schwarz again, we find
$$(c)~~~ \int_M \xi^2 u^4 \le 2S \big( \int_M u^2 |d\xi|\big)^2 + 8S (\int_M
u^2) \int_M \xi^2 |du|^2.$$

Plug (c) into (a) to get

$$\big( 1 - 8S \int_M u^2 \big) \int_M \xi^2 |du|^2 \le 2S \big(
\int_M u^2 |d\xi|\big)^2 + \int_M u^2 |d\xi|^2.$$

We now assume that $8S \int_M u^2 < 1$ and we choose a family of
functions $\xi_R$ such that $\xi_R$ is equal to $1$ in $B(x_0,R)$
(the ball with radius $R$ centered at some $x_0 \in M$), $\xi_R$ is
equal to $0$ outside $B(x_0,2R)$ and $|d\xi_R| \le 2/R$. Letting $R$
tend to infinity and using the fact that $\int_M u^2$ is finite, we
obtain that $du = 0$. Since $M$ has infinite volume, it follows that
$u=0$.  \qed

\subsection{Higher dimension, $M^n \looparrowright \HH^n \times \R, n \ge 3$}

Recall the formula for the Jacobi operator,
\begin{equation*}\label{E-dim3a-1}
J_M := - \Delta_M - \big( |A_M|^2 + \rich(N_M) \big)
\end{equation*}
where $N_M$ is a unit normal field along $M$ and $A_M$ the second
fundamental form of $M$ with respect to $N_M$
(Section~\ref{SS-jac}). \bigskip

Let $v_M := \hat{g}(N_M,\partial_t)$ be the vertical component of
the unit normal vector $N_M$. A simple computation gives that $\rich
(N_M) = - (n-1)(1-v_M^2)$. It follows that the Jacobi operator of
$M$ is given by
\begin{equation}\label{E-dim3a-2}
J_M := - \Delta_M + (n-1) (1-v_M^2) - |A_M|^2.
\end{equation}

We have the following theorem.

\begin{thm}\label{T-dim3}
Let $M^n \looparrowright \HH^n \times \R$ a complete oriented
minimal immersion. Assume that $M$ has \emph{finite total
curvature}, \ie $\int_M |A_M|^n \, d\mu_M < \infty.$
\begin{enumerate}
    \item For $n \ge 2$, the second fundamental form $A_M$ tends to
    zero uniformly at infinity.
    \item For $n \ge 3$, the Jacobi operator of the immersion has finite index and, more
    precisely, there exists a universal constant $C(n)$ such that
    \begin{equation}\label{E-dim3a-3}
    \mathrm{Ind}(J_M) \le C(n) \int_M |A_M|^n \, d\mu_M .
    \end{equation}
\end{enumerate}
\end{thm}

\textbf{Remarks.}\\[3pt]
(i) The examples $\cM_d$ prove that the converse statements in the
previous theorems are not true in general, see Section
\ref{SS-dim3-trans2}. \\
(ii) Note that we state the second assertion of Theorem~\ref{T-dim3}
only for $\mathrm{dim}(M) \ge 3$ (our proof does not apply in
dimension 2, see \cite{BB90}).\bigskip

\pf As in the proof of Theorem \ref{T-dim2}, the manifold $M$
satisfies a Sobolev inequality of the form (\ref{E-dim2-sob1}),
namely
$$\|f\|_{n/(n-1)} \le S(M) \|df\|_1 \text{ ~for all~ } f \in
C_0^{1}(M)$$

for some constant $S(M)$. Furthermore, the second fundamental form
$A_M$ satisfies the following Simons' equation (compare with
(\ref{E-dim2-sim3})),
$$-\Delta |A_M| \le |A_M|^3 + C(n) |A_M|,$$

for some constant $C(n)$ which only depends on the dimension (this
follows from the expression of the term $\Rh(A)$ as given in
\cite{Si68}. \bigskip

The de Giorgi-Moser-Nash technique applies (see \cite{BCS98},
Theorem 4.1) and it follows that $|A_M|$ tends to zero uniformly at
infinity.
\bigskip

Since $|v_M| \le 1$, the operator $J_M$ is bounded from below and
essentially self-adjoint. Furthermore, its index is less than or
equal to the index of the operator $- \Delta - |A_M|^2$ which is
also bounded from below and essentially self-adjoint. The estimate
(\ref{E-dim3a-3}) then follows by applying Theorem 39 in
\cite{BB90}. \qed \bigskip

\textbf{Remark.} The preceding results can be generalized to minimal
hypersurfaces in $\HH^n \times \R^k$ or $\HH^n \times \HH^k$.

\section{Applications}\label{S-appli}

In this section, we use the examples constructed in Section
\ref{S-examp} as barriers to prove some general results on minimal
hypersurfaces in $\HH^n \times \R$. These results generalize results
obtained in \cite{NSST08} for dimension $2$. Similar results hold
for $H$-hypersurfaces as well, see \cite{NSST08,BS08b}.\bigskip

\begin{thm}\label{T-appli-1min}
Let $\Gamma \subset \HH^n$ be a compact embedded hypersurface and
consider two copies of $\Gamma$ in different slices, $\Gamma_{-} =
\Gamma \times \{-a\}$ and $\Gamma_{+} = \Gamma \times \{a\} \subset
\HH^n \times \R$, for some $a > 0$. Assume that $\Gamma$ is convex.
\medskip

Let $M \subset \HH^n \times \R$ be a compact immersed minimal
hypersurface such that $\partial M = \Gamma_{-} \cup \Gamma_{+}$.
Then,
$$2a < \frac{\pi}{n-1} \text{ ~(the height of the family
of catenoids)}.$$ Furthermore, if $M$ is embedded,
\begin{enumerate}
    \item $M$ is symmetric with respect to the slice $\HH^n \times
    \{0\}$.
    \item The parts of $M$ above and below the slice of symmetry are
    vertical graphs.
    \item If $\Gamma$ is a horizontal graph and symmetric with respect to a hyperbolic
    hyperplane $P$, then $M$ is a horizontal graph and symmetric with respect to the
    vertical hyperplane $P \times \R$. In particular, if $\Gamma$ is an
    $(n-1)$-sphere, then $M$ is part of a catenoid.
\end{enumerate}
\end{thm}

\pf We reason ab absurdo. Suppose that the height of $M$ is greater
than or equal to $\displaystyle \frac{\pi}{n-1},$ that is $2a\ge
\displaystyle \frac{\pi}{n-1}$. We recall that the height of the
family of $n$-dimensional catenoids $\{\cC_{\rho},\; \rho\in (0,
\infty)\}$ is bounded from above by $\displaystyle \frac{\pi}{n-1},$
but each catenoid $\cC_\rho$ has height strictly less than
$\displaystyle \frac{\pi}{n-1}$. Now as $M$ is compact, there is a
(hyperbolic) radius $\rho_0$ big enough such that $M$ is strictly
contained inside the vertical cylinder  ${\mathbb M}_{\rho_0}$ of
radius $\rho_0$ (where ${\mathbb M}_{\rho_0}$ is a cylinder over a
$n-1$ sphere ${\mathbb S}_{\rho_0}\subset \HH^2\times\{0\}$ of
radius $\rho_0$) containing $M$ in its mean convex side. Recall
that, by the geometry of the catenoids, the catenoid $\cC_{\rho_0}$
whose distance to the $t$-axis is $\rho_0$ is contained in the
closure of the non mean convex side of ${\mathbb M}_{\rho_0}$
touching ${\mathbb M}_{\rho_0}$ just along the $n-1$ sphere
${\mathbb S}_{\rho_0}$. Hence, $M$ is  strictly contained in the
connected component of $\HH^2 \times \R \setminus \cC_{\rho_0}$ that
contains the $t$ axis of $\cC_{\rho_0}.$ Notice that the whole
family of catenoids $\cC_{\rho}$ is strictly contained in the slab
of $\HH^2\times \R$ with boundary $\Gamma_{-}\cup \Gamma_{+}$.
Starting from $\rho=\rho_0$, making $\rho\rightarrow 0,$ that is
moving the family of catenoids $\{\cC_{\rho},\, \rho\le \rho_0\}$
towards $M$, we will find a first interior point of contact with
some $\cC_{\rho}$ and $M$, since the family of catenoid cannot touch
the boundary of $M$. We arrive at a contradiction, by the the
maximum principle. The proof of the first part of the statement is
completed.\bigskip

Now using the family of  slices $\HH\times\{t\}\cup \HH\times\{-t\}$
coming from the infinity towards $M$ we get, by maximum principle,
that  $M$ is entirely contained in the closed slab whose boundary is
the slices $\HH\times\{a\}\cup \HH\times\{-a\}$ and
$\left(\HH\times\{a\}\cup \HH\times\{-a\}\right )\cap M=\partial
M.$\bigskip

In the same way, considering the family of vertical hyperplanes, we
get that $M$ is contained in the mean convex side of the vertical
cylinder ${\mathbb M}_\Gamma$ over $\Gamma$ and ${\mathbb
M}_\Gamma\cap M=\partial M.$ Now using Alexandrov Reflection
Principle on the slices, moving the slices from $t=a$ towards $t=0$,
by vertical reflections, we get that the reflection of $M^+=M\cap\{
t\geqslant 0\}$, the part of $M$ above $t=0$, on the horizontal
slice $t=0$,  is above $M^-=M\cap \{t\leqslant 0\}$. Moreover, we
find that $M^+$ is a vertical graph. In the same way, moving the
slices from $t=-a$ towards $t=0$, doing vertical reflections, we get
that the reflection of $M^-=M\cap \{t\leqslant 0\}$, the part of $M$
above $t=0$, on the horizontal slice $t=0$, is below $M^+=M\cap
\{t\ge 0\}$. We conclude that $M^-=M^+$, hence both $M^+$ and $M^-$
are vertical graphs and $M$ is symmetric with respect to the slice
$\HH^n \times \{0\}$.  Therefore, the proof of the second part of
the Statement is  completed.\bigskip

Let us assume now that $P\subset\HH\times\{0\}$ is a hyperplane of
symmetry of $\Gamma.$  Consider the vertical hyperplane ${\mathbb
P}=P\times \R$ and the family of hyperplanes ${\mathbb P}_t$ at
signed distance $t$ from ${\mathbb P}$ obtained from ${\mathbb P}$
by horizontal translations along an oriented geodesic $\gamma$
orthogonal to ${\mathbb P}$ at the origin. Choosing $|t|$ big
enough, we move the family ${\mathbb P}_t$ towards ${\mathbb P}$ (in
the two sides of $\HH^2 \times \R \setminus {\mathbb P}$), doing
Alexandrov Reflection Principle on ${\mathbb P}_t$, taking into
account that $\Gamma$ is a horizontal graph and that the symmetric
of $\partial M$ on ${\mathbb P}_t$ stays on the slices $t=\pm a$, so
that it does not touch the interior of $M$. We can argue as before
to conclude that ${\mathbb P}$ is a hyperplane of symmetry of $M$.
Of course, is $\Gamma$ is rotationally symmetry then $M$ is a
minimal hypersurface of revolution. Henceforth, by the
classification theorem, $M$ is part of a catenoid. This completes
the proof of the theorem. \qed


\bigskip

\begin{thebibliography}{10}

\bibitem{BGS87}
Lucas Barbosa, Jonas Gomes, and Alexandre Silveira.
\newblock Foliation of $3$-dimensional space forms by surfaces with constant
  mean curvature.
\newblock {\em Bol. Soc. Bras. Mat.}, 18:1--12, 1987.

\bibitem{BB90}
Pierre B\'{e}rard and G\'{e}rard Besson.
\newblock Number of bound states and estimates on some geometric invariants.
\newblock {\em J. Funct. Anal.}, 94:375--396, 1990.

\bibitem{BCS97}
Pierre B\'{e}rard, Manfredo do~Carmo, and Walcy Santos.
\newblock The index of constant mean curvature surfaces in hyperbolic
  {$3$}-space.
\newblock {\em Math. Z.}, 224:313--326, 1997.

\bibitem{BCS98}
Pierre B\'{e}rard, Manfredo do~Carmo, and Walcy Santos.
\newblock Complete hypersurfaces with constant mean curvature and finite total
  curvature.
\newblock {\em Ann. Global Anal. Geom.}, 16:273--290, 1998.

\bibitem{BSA09L}
Pierre B\'{e}rard and Ricardo Sa~Earp.
\newblock Lindeloef's theorem for catenoids, revisited.
\newblock {\em {arXiv:0907.4294}}, 2009.

\bibitem{BS08b}
Pierre B\'{e}rard and Ricardo {Sa Earp}.
\newblock {$H$-hypersurfaces in {$\HH^n \times \R$} and applications}.
\newblock {\em Mat. Contemp.}, 34, 2008.

\bibitem{CP80}
Manfredo~do Carmo and C.K. Peng.
\newblock Stable complete minimal hypersurfaces.
\newblock In {\em Proceedings Beijing Symposium on Diff. Geom. and Diff. Eq.,
  ed. by S.S. Chern and W.W. Tsun}, pages 1349--1358, 1980.

\bibitem{CM02}
Tobias~H. Colding and William~P. Minicozzi.
\newblock Estimates for parametric elliptic integrands.
\newblock {\em IMRN International Mathematics Research Notices}, 6:291--297,
  2002.

\bibitem{CR07}
Pascal Collin and Harold Rosenberg.
\newblock Construction of harmonic diffeomorphisms and minimal graphs.
\newblock {\em {arXiv:0701.547v1}}, 2007.

\bibitem{LR98}
Levi~Lopes de~Lima and Wayne Rossman.
\newblock On the index of constant mean curvature {$1$} surfaces in hyperbolic
  space.
\newblock {\em Indiana Univ. Math. J.}, 47:685--723, 1998.

\bibitem{CS90}
Manfredo~P. do~Carmo and Alexandre~M. Da~Silveira.
\newblock Index and total curvature of surfaces with constant mean curvature.
\newblock {\em Proc. Amer. Math. Soc.}, 110:1009--1015, 1990.

\bibitem{FC82}
Doris Fischer-Colbrie.
\newblock On complete minimal surfaces with finite {M}orse index in
  three-manifolds.
\newblock {\em Invent. Math.}, 82:121--132, 1985.

\bibitem{FCS80}
Doris Fischer-Colbrie and Richard Schoen.
\newblock The structure of complete stable minimal surfaces in {$3$}-manifolds
  of nonnegative scalar curvature.
\newblock {\em Comm. Pure Appl. Math.}, 33:199--211, 1980.

\bibitem{HR08}
Laurent Hauswirth and Harold Rosenberg.
\newblock Minimal surfaces of finite total curvature in {$\HH \times \R$}.
\newblock {\em Mat. Contemp.}, 31:65--80, 2006.

\bibitem{HS74}
David Hoffman and Joel Spruck.
\newblock Sobolev and isoperimetric inequalities for riemannian submanifolds.
\newblock {\em Comm. Pure Appl. Math.}, 27:715--727, 1974.

\bibitem{HS75}
David Hoffman and Joel Spruck.
\newblock A correction to : Sobolev and isoperimetric inequalities for
  riemannian submanifolds.
\newblock {\em Comm. Pure Appl. Math.}, 28:765--766, 1975.

\bibitem{Law80}
H.~Blaine Lawson, Jr.
\newblock {\em Lectures on minimal submanifolds. {V}ol. {I}}, volume~9 of {\em
  Mathematics Lecture Series}.
\newblock Publish or Perish Inc., Wilmington, Del., second edition, 1980.

\bibitem{Lin870}
Lorenz Lindeloef.
\newblock Sur les limites entre lesquelles le cat\'{e}no\"{\i}de est une
  surface minimale.
\newblock {\em Math. Annalen}, 2:160--166, 1870.

\bibitem{NSST08}
Barbara Nelli, Ricardo Sa~Earp, Walcy Santos, and Eric Toubiana.
\newblock Uniqueness of {H}-surfaces in {$\HH^2 \times \R, |H| \le 1/2$}, with
  boundary one or two parallel horizontal circles.
\newblock {\em Ann. Global Anal. Geom.}, 33:307--321, 2008.

\bibitem{Sa08}
Ricardo Sa~Earp.
\newblock Parabolic and hyperbolic screw motion surfaces in {$\HH\sp
  2\times\R$}.
\newblock {\em Journal of the Australian Math. Soc.}, 85:113--143, 2008.

\bibitem{ST05}
Ricardo Sa~Earp and Eric Toubiana.
\newblock Screw motion surfaces in {$\HH\sp 2\times\R$} and {$\Ss\sp
  2\times\R$}.
\newblock {\em Illinois J. Math.}, 49:1323--1362, 2005.

\bibitem{ST08}
Ricardo Sa~Earp and Eric Toubiana.
\newblock An asymptotic theorem for minimal surfaces and existence results for
  minimal graphs in {$\HH\sp 2\times\R$}.
\newblock {\em Math. Annalen}, 342:309--331, 2008.

\bibitem{ST09}
Ricardo Sa~Earp and Eric Toubiana.
\newblock Minimal graphs in {$\HH\sp n\times \R$} and {$\R^{n+1}$}.
\newblock {\em arXiv:0908.4170}, 2009.

\bibitem{SSY75}
Richard Schoen, Leon Simon, and Shing-Tung Yau.
\newblock Curvature estimates for minimal hypersurfaces.
\newblock {\em Acta Math.}, 134:275--288, 1975.

\bibitem{Si68}
James Simons.
\newblock Minimal varieties in riemannian manifolds.
\newblock {\em Ann. of Math.}, 88:62--105, 1968.

\bibitem{V50}
Georges Valiron.
\newblock {\em {\'{E}quations fonctionnelles. Applications}}.
\newblock Masson, 1950.

\end{thebibliography}

\vspace{1cm}

\begin{tabular}{lll}
Pierre B\'{e}rard & \hphantom{xxxxxxxxxxxxxx}& Ricardo Sa Earp\\
Universit\'{e} Joseph Fourier && Pontif{\'{i}}cia Universidade Cat\'{o}lica\\
Institut Fourier (\textsc{ujf-cnrs}) && do Rio de Janeiro \\
B.P. 74 && Departamento de
Matem\'{a}tica\\
38402 Saint Martin d'H\`{e}res Cedex && 22453-900 Rio de Janeiro - RJ\\
France && Brazil\\
\verb+Pierre.Berard@ujf-grenoble.fr+ && \verb+earp@mat.puc-rio.br+\\
\end{tabular}\bigskip

\end{document}